\documentclass{article}
\usepackage{arxiv}
\usepackage[utf8]{inputenc}
\usepackage{amsmath}
\usepackage{amssymb,amsfonts,textcomp}
\usepackage[T1]{fontenc}
\usepackage[english]{babel}
\usepackage{xcolor}
\usepackage{color}
\definecolor{hellgrau}{RGB}{223,223,223}
\definecolor{hellgelb}{RGB}{255,255,199}
\definecolor{helllila}{RGB}{209,175,239}
\definecolor{hellblau}{RGB}{210,240,242}
\definecolor{midblau}{RGB}{180,230,232}
\definecolor{midgrau}{RGB}{174,174,174}
\usepackage{colortbl}\definecolor{hellgrau}{RGB}{223,223,223}
\usepackage{array}
\usepackage{booktabs}       
\usepackage{url}
\usepackage{hhline}
\usepackage{hyperref}
\usepackage{nicefrac}
\usepackage{graphicx}
\usepackage{float}
\usepackage{lscape}

\usepackage[backend=biber,
style=numeric,
sorting=nyt,
giveninits=true]{biblatex}
\title{On the evidence of a trigonometric function value system in Babylon}
\author{Jens Kleb \\ }
\date{22 December 2021}

\begin{document}

\maketitle

\begin{abstract}
Ratios and coefficients are used to simplify calculations. For geometric usage these values also called function values. Like in Egypt also in Babylon such a value system can be shown.

The reconstructed calculation sequence, of the Plimpton 322 cuneiform tablet, presented and described here, shows in its completeness that, around 3800 years ago there already was a systematically applied exact measuring system, with the usage of trigonometric function values.

With this approach one can plausibly explain that, as we still it practice today, a geometry of the circle has been used for this calculation. It is based on, but not only, of the usage of the regularities, which 1200 years later, were named after Pythagoras. 

During a second calculation step, for an intended scaled documentation, presentation or transfer to other locations, the dimensionless calculated function value, was extended, with a fractional part, or more exact spelled a common factor. This transformation creates a real measurable length from the ratios, always related to the smallest unit.

The systematic usage by means of inclining triangles already in old Babylonian times, goes far beyond the previously known level of this age.

The accuracy of the individual calculation is just as verifiable and on exactly the same level as the trigonometric functions used today. However, at least at this point in time, the Babylonians were content with dividing a quarter circle into at least 150 inclining triangles.

This old Babylonian trigonometric function value system, together with the here described
\begin{center}
	\textit{"Babylonian diagonal calculation"},
\end{center}
is thus the forerunner of the Greek chord calculation and the Indo-Arabic trigonometric functions of the early Middle Ages, which we still use today. 

In addition, the new approach to Plimpton 322 plausible explains why and how such values could have been used.

\end{abstract}

\newpage
\tableofcontents

\newpage
\section{Function values in Egypt and Babylon, a short Overview}
\subsection{Proportions and Function Values, Evidence and Usage}
Function values and their use are not uncommon. Even millennia ago they were used, as ratio values to describe weights per unit of volume, or also inclination ratios to simplify practical applications. The Rhind papyrus shows us this for Ancient Egypt \cite{chase1929}. On this papyrus is a description of the "Seked", a coefficient that corresponds to our cotangent today. Mathematical calculations on cuneiform tablets, from the old Babylonian period confirm that such ratio values and coefficients were also used in Mesopotamia \cite{robson2014}(BM085194; BM96957; BM085200; BM96954 and others). A trigonometric function value reflects a dimensionless ratio between two sides of a triangle. The usage of a ratio, simplifies the erection of structures, or the description of slopes using the value related to a standard dimension of length.
As already mentioned by D. Fowler \cite{fowler1998}, these are also known as coefficients and have been used for a wide variety of calculations. The approximation mentioned there, at \textit{YBC7289} for the value of the $\sqrt{2}$ is also explicitly referred as a coefficient.

Whether as a coefficient, ratio, or ultimately as a function value, they all describe the same thing. 
Therefore, the coefficient $\sqrt{2}$ is at the same time, also the function value known today as Secant $45$ Degree. 
In contrast to the ancient Egyptian Seked as today's functional value of the cotangent\cite{Maor2002}, the geometric concept of the ``diagonal'' was often used in the ancient Babylonian region, corresponding to today's Secant \cite{Joseph2000}.
For use in immediate calculation or documentations, this ratio has to be multiplied by the required, or with any other practicable, length. 

In later times, fixed numbers were assigned to the ratio or function values, which we know today as ``angles''\cite{bowen2020}.

While the use of gradient ratios and thus a functional value at the latest since the Rhind papyrus\cite{chase1929} has been accepted for Egypt, this is currently not the case for Mesopotamia and the almost equally dated Old Babylonian period.
\\

\label{img001.jpg}
\begin{figure}[h]
	\centering
	\includegraphics{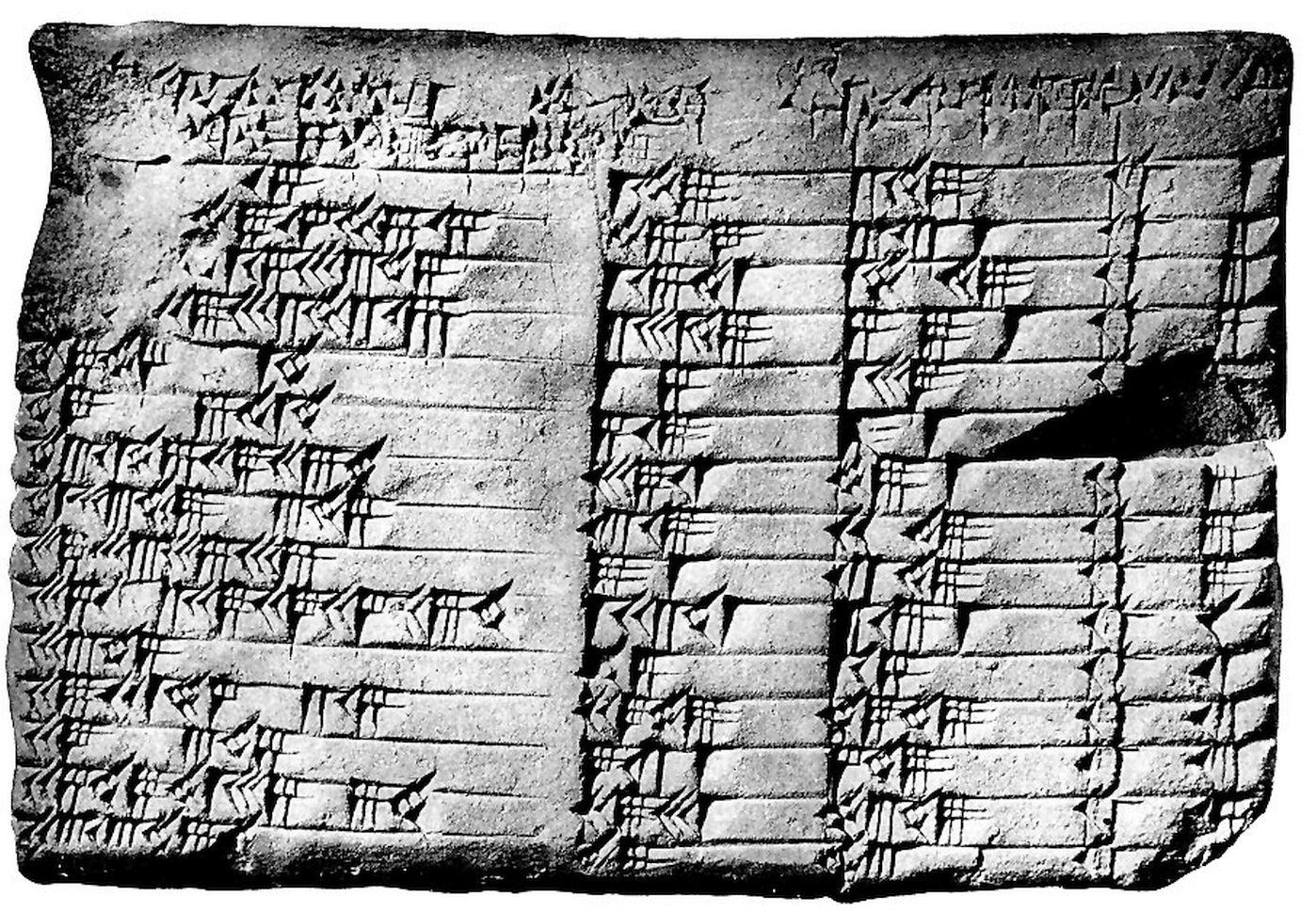}
	\caption{Photography of Plimpton 322 from Neugebauer and Sachs \cite{neugebauer1945}.}
\end{figure}

\subsection{Previous research on Plimpton 322}

This particular cuneiform tablet has been part of the Plimpton Collection of the Columbia University New York since 1923. It was probably found, or bought, by E. J. Banks at a location near Tell Senkereh the ancient Larsa. The age of the clay tablet is estimated at 3800 years. The fact that, the cuneiform numbers on the obverse side form right-angled triangles in each of its 15 lines, was first recognized by O. Neugebauer and A. Sachs \cite{neugebauer1945}. The so-called Pythagorean triples\cite[36]{Herrmann2014} and the implied use of the Pythagorean equation $a^2 + b^2 = c^2$ has attracted great interest in science since then.

Many theories about the calculation and use of the cuneiform table have been discussed intensively, but some questions always remain unanswered. Most recently, Buck \cite{buck1980}, Friberg \cite{friberg1981} and \cite{friberg2007}, Høyrup \cite{hoyrup1999}, Robson \cite{robson2001}, Britton \cite{britton2011}, as well as Mansfield and Wildberger \cite{mansfield2017} and \cite{mansfield2021}, have dealt with Plimpton 322.

Two opposing scientific parties emerged. While Robson \cite{robson2001} and Friberg \cite{friberg2007} suggest the goal of Plimpton 322 is in the generation of Pythagorean triples (also\cite{bruins1949},\cite{bruins1967},\cite{schmidt1980},\cite{Muroi2013}) specifically it as ``an approach to the creation of reciprocal pairs''.
Eleanor Robson is known as one of the most renowned scientists in the field of mathematical cuneiform texts and so she
has defined \cite[p. 176]{robson2001} six criteria for an interpretation of Plimpton 322, which need to be met to be an acceptable reconstruction. The six criteria were formulated to be used against the interpretation of angles or incline proportions in ancient Babylonian times. They will also be examined in this paper \{\ref{section:criteria}\}. It will be shown that, with the suggested calculation sequence and its results described in this paper, there is absolutely no contradiction with them. The content of the cuneiform tablet can completely fulfill all of these criteria.
On the other hand, Joseph\cite{Joseph2000}, Hajossy\cite{Hajossy2016} and finally Mansfield and Wildberger \cite{mansfield2017} believe that it was ``a kind of trigonometric table'' which might have served as the basis for further interpolations. In their paper, Mansfield and Wildberger conclude, that no other previous explanation has such a high degree of plausibility. 
Because of that, we have to look at this last mentioned paper\cite{mansfield2017}, in more Detail:

Unfortunately, it remains open, however, why and how the reconstructed values $b$ and $d$ contained in table 13 of Mansfield and Wildberger \cite{mansfield2017} have been calculated, which are the only ones we can find in real on the Plimpton 322 cuneiform tablet. Especially, as the values $\beta$ and $\delta$ of the same table 13 already define the inclination / angle exactly.

According to Mansfield and Wildberger the reconstruction of calculation for contents of the table relates, according Mansfield and Wildberger, to a sequence of ``11 steps'' \cite[pp. 405, 406]{mansfield2017}. They assume that there was a maximum number of 38 resulting lines and that it was used as a trigonometric table. However further questions arise. According to the given procedure by the authors, the values must be ``sorted'' in descending order in the 7th step. Why descending and not ascending, as in our currently used trigonometric tables, is not explained. The content of Plimpton 322 by this presented procedure, should have been therefore only a compilation, created in the last ``step 11'' from all interim results. According to them it was not a step-by-step calculation procedure. Why for such a longer lasting compilation of data, are the known errors acceptable, is the remaining question.

However, one of the most crucial points of the Mansfield and Wildberger reconstruction is step 8. This step does not give a plausible explanation for the necessity of the values given in column I of the cuneiform tablet and, to make matters worse, it is added in the explanation of step 8, that the authors do not believe that the values $\beta$ and $\delta$ were calculated from the given column I.

In the given procedure, step 8 is just somewhat like a necessary insertion, to fulfill the given columns of the tablet. It is just there to offer the values of column I, which is otherwise not imminently necessary for their procedure. It is evident that, from the values in step 5, which should be sorted again in step 7, during step 9 the values of column II for value $b$ and column III for value $d$, should be calculated. At the end, nevertheless, the question is now, column I would be irrelevant for the whole calculation and why should this be correct?

In conclusion, it can be stated, that the calculation principles of Plimpton 322 which have been sought since Neugebauer and Sachs 1945 \cite{neugebauer1945} and also been requested by Scriba and Schreiber \cite{scriba2001}, \cite{scriba2015}, have not yet been satisfactorily described\cite[4]{Phillips2011}.
\, \

\section{Babylonian Calculus}
\label{section:calculus}

The Babylonian calculus is quite different from the modern. Numbers are represented in sexagesimal form, i.e. to the base of 60. The places in transcribed numbers are separated by full stops. The decimal point and trailing places of zeroes are unknown. This means that a Babylonian sexagesimal number $s$ is not a single number, it merely represents in fact a whole set of numbers. The set consists of all numbers $s \cdot 60^n$ for all integers $n$. So $1$ is also $60$, $3600$, $216000$ and so forth. If one wants to add or subtract $1$, one has to realize from the context which representative has to be added.(for further explanations see:\cite{Lemmermeyer2015}\cite{Lemmermeyer2015a}\cite{lehmann1992})

\newpage
\section{The trigonometric calculation with diagonals in old Babylonian times (in general)}
\label{section:definitions}
\begin{figure}[h]
	\centering
	\includegraphics[scale=0.2]{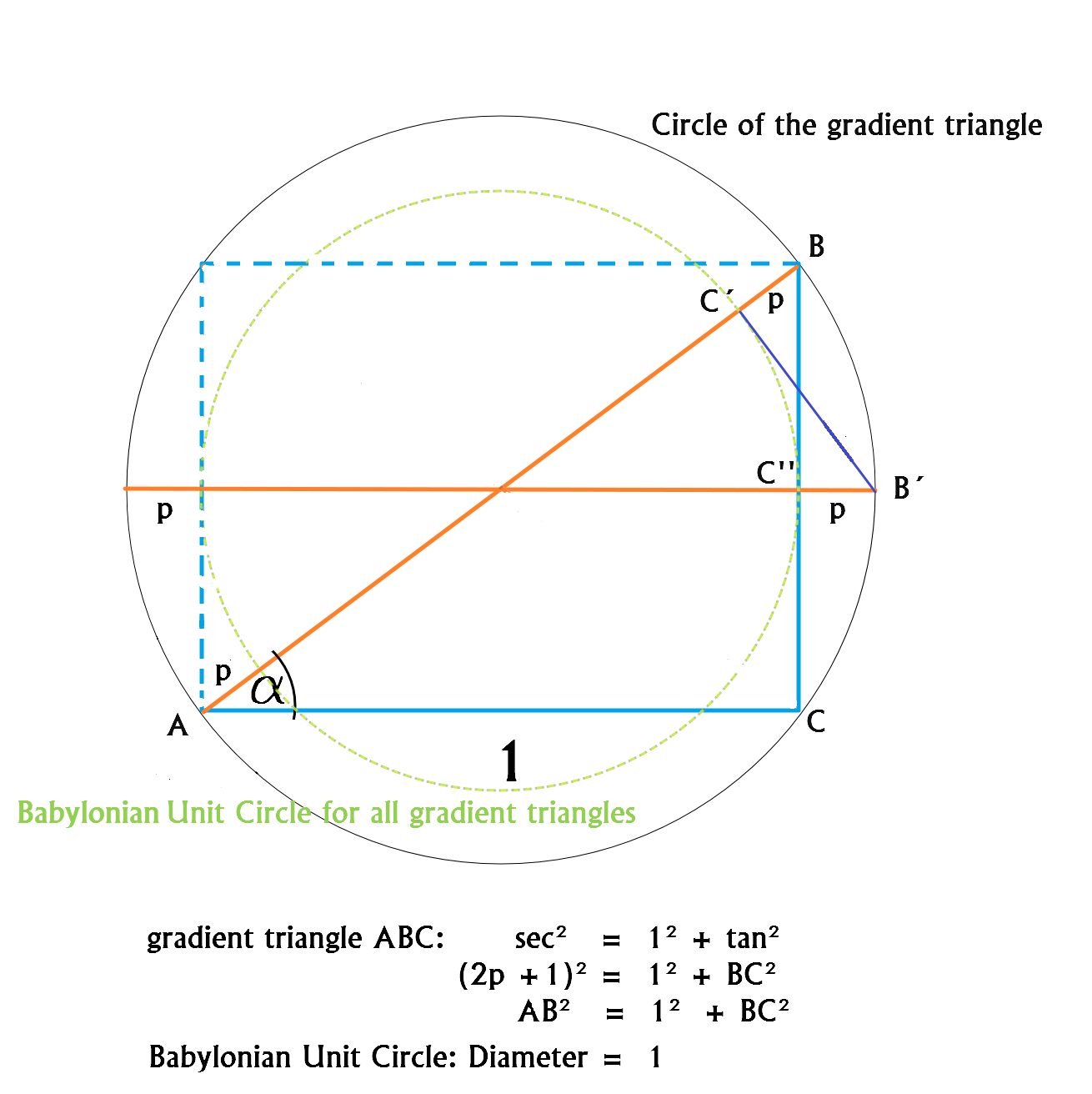}
	\caption{Geometric construct for the calculation sequence.}
\end{figure}

\begin{table}[h!]
	\centering
	\caption{Definitions at the triangle and its connected circles}
	\begin{tabular}{p{0.15\textwidth} p{0.1\textwidth}p{0.1\textwidth}p{0.15\textwidth}p{0.1\textwidth}p{0.1\textwidth}p{0.1\textwidth}}
		\toprule
		\centering\textbf{Side name in the triangle} & \centering \textbf{Babylonian name} & \centering \textbf{Side name} & \centering \textbf{Section between points} & \centering \textbf{Proportion} & \centering \textbf{Proportion in Bab. Calculation} & \centering \textbf{Trig. Function}\tabularnewline
		\midrule
		\centering Hypotenuse & \centering Diagonal & \centering $c$ & \centering $\overline{AB}$ & \centering $\dfrac{\textrm{Hypotenuse}}{\textrm{Adjacent}}$ & \centering $\dfrac{c}{b}$ & \centering Secant\tabularnewline\tabularnewline
		\centering Opposite & \centering Width &  \centering $a$ & \centering $\overline{BC}$ & \centering $\dfrac{\textrm{Opposite}}{\textrm{Adjacent}}$ & \centering $\dfrac{a}{b}$ & \centering Tangent\tabularnewline\tabularnewline
		\centering Adjacent & \centering Length & \centering $b$ & \centering $\overline{AC}$ is fixed to sexagesimal 1 & \centering $\dfrac{\textrm{Adjacent}}{\textrm{Hypotenuse}}$ & \centering $\dfrac{b}{c}$ & \centering \tabularnewline\tabularnewline
		\centering Height of the chord to its arc segment & \centering Arrow & \centering $p$ & \centering $\overline{C''B'}=\overline{C'B}$ & \multicolumn{2}{ p{0.2\textwidth} }{ \raggedright Section outside the Bab. Unit Circle to the Gradient Circle (unique
			for each gradient triangle and implicit to its enclosed angle $\alpha$).} & \centering $\dfrac{\textrm{Exsecant}}{\textrm{2}}$ \tabularnewline\midrule
		\multicolumn{7}{p{0.95\textwidth}}{To distinguish the modern function value from the value multiplied
			by 60, we use the term \textbf{\textit{bab-sec}} and \textbf{\textit{bab-tan}}. For the Babylonian calculation
			the Adjacent $b$ is fixed to sexagesimal 1, which is 60 as decimal
			value.}\cite{Gellert1969}\cite{Bronstein1964}\\\bottomrule
	\end{tabular}
\end{table}

The length $p$ is in a relationship to the hypotenuse of the triangle, which is equal to diagonal of the inscribed chordal quadrilateral (created by doubling the triangle). The diagonal is also the diameter of the circumcircle. This diameter minus triangle side $b$ and divided by 2 is the Length of $p$ the arrow. Because $\overline{AC}$ is fixed to $1$ for this construct, it is
\[
\text{therefore  },p = (\text{sec}\,\alpha - 1 ) / 2  \  \  \ \text{  ;  }\  \ \  \  \ \text{  if } \overline{AC} = b = 1 ;   
\]
The diagonal of the gradient triangle and at the same time the function value $secant$ of the included angle $\alpha$ is:
\[
\text{diagonal} = \text{sec}\,\alpha = 1 + 2 p \  \  \ \text{  ;  } 2p = exsec  \ \alpha
\]
The fixed side length (adjacent cathetus, or named as "base")by $1$, is strictly defined for any inclining/gradient triangle.

This also results in other relationships, some of which are unused or completely forgotten today for the unit circle, for instance:
\[
1 = \text{cos}\, \alpha + 2 p \times \text{cos}\, \alpha
\]
According to Thales's theorem, it is immediately noticeable that, , a right-angled triangle is created, with the adjacent cathetus $b = 1$, and the opposite cathetus $a = \text{tan}\, \alpha$,  above $c$ = the diagonal = secant $\alpha$. The exact proportions hereby depend on the gradient ratio, or our modern angle $\alpha$.

As Babylonian normality, this is done efficiently, completely correct and without squaring or taking square roots. Squaring was only used when determining the width (opposite or side $a$), as in column I of the cuneiform tablet Plimpton 322. It goes without saying that, the sine is also usable too for this, but this function value was not necessary. Because the side $a$ within the unit circle is equal to the function value of $tangent$ $\alpha$ it is valid: 
\[
sec^{2}\alpha -1^{2} = tan^{2}\alpha
\]

The Babylonians virtually ``shrinked'' the circle of the inscribed gradient triangle to a diameter of $1$ or $60$ for the Babylonian unit circle.

Always related to the most effective use, the secant ratio as the diagonal and the tangent ratio as the width or opposite side $a$, of the angle $\alpha$ at Point $A$, were certainly the preferred function values of the Babylonians for the documentation of gradient ratios.

We today use the tangent function in everyday life more often than all other function values, because it is easier to imagine. Whether specifying the gradient ratios of a staircase with an exemplary ratio of $1/3$, or specifying the gradient on street signs in percent, or the falling gradient of a river in per million, all are given in such proportions that correspond to the \textit{tangent function value} of the angle $\alpha$. 
That's why, this function value is defined by the ratio of width (here also named opposite side) $a$, divided by the length of the base side $b$.

The Babylonians viewed their unit circle as a whole (Babylonian Unit circle) and not, as we do today, essentially as a quadrant wise geometric solution based in the center of the circle. However, the results are the same because all proportions or gradient ratios of the sides to each other, remain also the same.

Based on this clear geometric construct established and used by the Babylonians, we see the already mentioned Pythagorean laws for the right-angled triangle and its application to the trigonometric functions within the Babylonian unit circle. Likewise, the Thales theorem is implicitly used, since every angle formed over the diameter of the circle must be right-angled one. Moreover, we recognize the regularities of chord calculation, because the side length $a$ of the opposite peripheral angle $\alpha$, is equal to that side length of the doubled centric-angle $\alpha$ at the circle center and its triangle.

Therefore, as is well known today, the following applies:
\begin{gather*}
	\text{sin} (2 \alpha) = 2\, \text{sin}\,\alpha \times \text{cos}\,\alpha
\end{gather*}

\begin{table}[H]
	\centering
	\caption{Reconstructed tablet values in sexagesimal numbers with corrections}
	\begin{tabular}{p{0.5cm}p{3cm}p{4.8cm}p{2.0cm}p{2.0cm}p{1.7cm}}
		\toprule
		&
		\centering \includegraphics[scale=0.9]{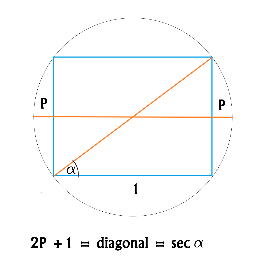}
		&
		\centering \includegraphics[scale=0.7]{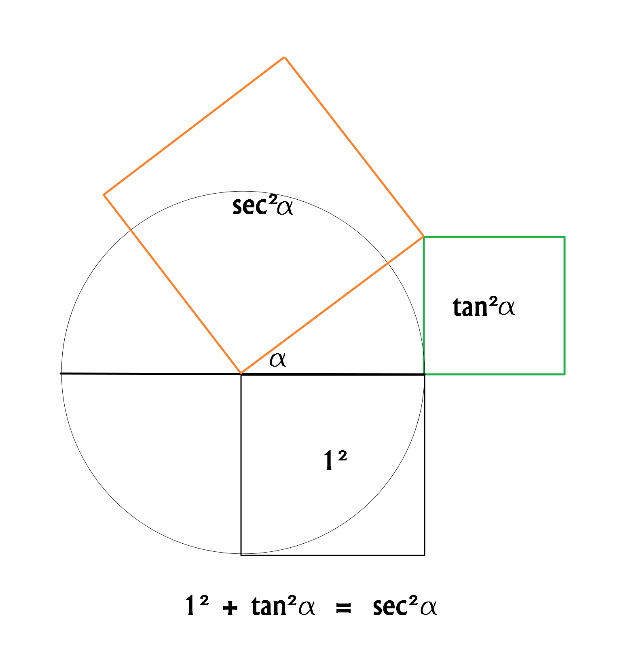}
		&
		\multicolumn{2}{p{4cm}}{\centering \includegraphics[scale=0.9]{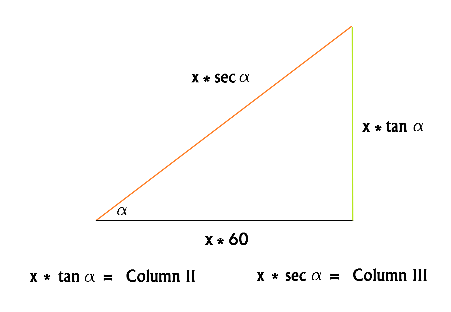}} 
		&
		\\\midrule
		\cellcolor{midgrau}
		&
		\cellcolor{hellgrau}
		\textbf{Column 0} (reconstr.)
		&
		\textbf{Column I}
		&
		\textbf{Column II}
		&
		\textbf{Column III}
		&
		\textbf{Column IV}
		\\
		\cellcolor{midgrau}
		&
		\cellcolor{hellgrau}
		Note the doubled value, later add 1 and square this sum
		&
		The intermediate (takiltum) is the square of the diagonal, which minus 1 becomes the square of the width
		&
		The Result for the width
		&
		The Result for the diagonal
		&
		It's Place
		\\\midrule
		\cellcolor{midgrau}
		&
		\cellcolor{hellgrau}
		24.30&
		1.59.00.15&
		1.59&
		2.49&
		1
		\\
		\cellcolor{midgrau}
		&
		\cellcolor{hellgrau}
		23.46.2.30&
		1.56.56.58.14.50.06.15&
		56.07&
		1.20.25\textsuperscript{1}&
		2
		\\
		\cellcolor{midgrau}
		&
		\cellcolor{hellgrau}
		23.6.45&
		1.55.07.41.15.33.45&
		1.16.41 &
		1.50.49 &
		3
		\\
		\cellcolor{midgrau}
		&
		\cellcolor{hellgrau}
		22.24.16 &
		1.53.10.29.32.52.16 &
		3.31.49 &
		5.9.1 &
		4\\
		\cellcolor{midgrau}
		&
		\cellcolor{hellgrau}
		20.50&
		1.48.54.01.40&
		1.05 &
		1.37 &
		5\\
		\cellcolor{midgrau}
		&
		\cellcolor{hellgrau}
		20.10&
		1.47.06.41.40 &
		5.19 &
		8.01 &
		6\\
		\cellcolor{midgrau}
		&
		\cellcolor{hellgrau}
		18.41.20&
		1.43.11.56.28.26.40 &
		38.11 &
		59.1 &
		7\\
		\cellcolor{midgrau}
		&
		\cellcolor{hellgrau}
		18.3.45&
		1.41.33.45.14.03.45\textsuperscript{1} &
		13.19 &
		20.49 &
		8\\
		\cellcolor{midgrau}
		&
		\cellcolor{hellgrau}
		16.54 &
		1.38.33.36.36 &
		8.01\textsuperscript{1}&
		12.49 &
		9\\
		\cellcolor{midgrau}
		&
		\cellcolor{hellgrau}
		15.33.53.20 &
		1.35.10.02.28.27.24.26.40 &
		1.22.41 &
		2.16.1&
		10\\
		\cellcolor{midgrau}
		&
		\cellcolor{hellgrau}
		15&
		1.33.45&
		45 &
		1.15&
		11\\
		\cellcolor{midgrau}
		&
		\cellcolor{hellgrau}
		13.13.30&
		1.29.21.54.02.15 &
		27.59&
		48.49 &
		12\\
		\cellcolor{midgrau}
		&
		\cellcolor{hellgrau}
		12.15&
		1.27.00.03.45 &
		2.41\textsuperscript{1} &
		4.49 &
		13\\
		\cellcolor{midgrau}
		&
		\cellcolor{hellgrau}
		11.45.20 &
		1.25.48.51.35.06.40 &
		29.31&
		53.49 &
		14\\
		\cellcolor{midgrau}
		&
		\cellcolor{hellgrau}
		10.40 &
		1.23.13.46.40&
		56 &
		1.46\textsuperscript{1} &
		15\\ \midrule
		\multicolumn{5}{m{14cm}}{
		\textsuperscript{1} corrected number}
		\\\bottomrule
	\end{tabular}
\end{table}

Without going into detail, it can be seen how trigonometric function systems developed during the following millennia.
 
This was plausibly initiated after cutting out the Arrow length $\textbf{\textit{p}}$, which has on the other hand initially created the subdivision of the geometric segments, from the unit circle. The development goes from the initial created trigonometric diagonal calculation by the Babylonians, of constantly reducing the reference triangles, first to the Greek chord calculation with double angles and then by halving again to our modern trigonometry.
 
However, in contrast to nearly all ancient weight or length systems, the main values of all these trigonometric function systems are as usable and accurate as the modern calculated ones.

\newpage
\section{The calculation procedure on Plimpton 322}
\subsection{Introduction}

By using the following description and the attached tables, every step of the calculation of the values, existent on the Plimpton 322 tablet, will be transparent and comprehensible.
 
Fundamentally, however, it must be noted in prior, that all our previous knowledge, of the science during this age is extremely incomplete. We should not doubt of the high level in ancient science and therefore new findings should be analyzed objectively at any time, especially if these due to the precise mathematics, leaving little room for speculation, as it is normal for this outstanding science at all times.

Moreover and even if this is not part of the traditional way of thinking up to now, we have to accept that in addition to the complex geometric imagination, at the time Plimpton 322 was written, not only the use of the Pythagorean theorem was known. Rather, and this underlines, not only the heading above column I, a complex function value system has already been used at this time, which will also be shown here \cite{Kleb2021b}\cite{Kleb2021a}.

The Evidence for this, is in addition given to us directly on the heading of column I, on the tablet via the properties of trigonometric Pythagoras\cite{Gellert1969}\cite{Bronstein1964}.
\\
\\

\begin{table}[H]
	\centering
	\caption{Reconstructed tablet values in decimal numbers.}
	\begin{tabular}{p{0.5cm}p{3cm}p{4.8cm}p{2.0cm}p{2.0cm}p{1.7cm}}
		\toprule
		\cellcolor{midgrau}&
		\cellcolor{hellgrau}
		\textbf{Column 0} (reconstr.)
		&
		\textbf{Column I}
		&
		\textbf{Column II}
		&
		\textbf{Column III}
		&
		\textbf{Column IV}
		\\
		\cellcolor{midgrau}
		&\cellcolor{hellgrau}
		Note the doubled value, later add 1\textsuperscript{\$} and square this sum
		&
		The intermediate (takiltum) is the square of the diagonal, which minus 1\textsuperscript{\$} becomes the square of the width
		&
		The Result for the width
		&
		The Result for the diagonal
		&
		It´s Place
		\\\midrule
		\cellcolor{midgrau}
		&\cellcolor{hellgrau}
		24,5 &
		1,983402$p$7 &
		119 &
		169 &
		1 \\
		\cellcolor{midgrau}
		&\cellcolor{hellgrau}
		23,76736$p$1 &
		1,949158552088* &
		[16835]
		
		3367 &
		[24125]
		
		4825\textsuperscript{1} &
		2\\
		\cellcolor{midgrau}
		&\cellcolor{hellgrau}
		23,1125 &
		1,918802126736$p$1 &
		4601 &
		6649 &
		3\\
		\cellcolor{midgrau}
		&\cellcolor{hellgrau}
		22,40$p$4 &
		1,886247906721* &
		12709 &
		18541 &
		4\\
		\cellcolor{midgrau}
		&\cellcolor{hellgrau}
		20,8$p$3 &
		1,815007716049* &
		[325]
		
		65&
		[485]
		
		97 &
		5\\
		\cellcolor{midgrau}
		&\cellcolor{hellgrau}
		20,1$p$6 &
		1,785192901234* &
		319 &
		481 &
		6\\
		\cellcolor{midgrau}
		&\cellcolor{hellgrau}
		18,6$p$8 &
		1,719983676268* &
		2291 &
		3541 &
		7\\
		\cellcolor{midgrau}
		&\cellcolor{hellgrau}
		18,0625 &
		1,692709418402$p$7\textsuperscript{1} &
		799 &
		1249 &
		8\\
		\cellcolor{midgrau}
		&\cellcolor{hellgrau}
		16,9 &
		1,642669$p$4 &
		481\textsuperscript{1} &
		769 &
		9\\
		\cellcolor{midgrau}
		&\cellcolor{hellgrau}
		15,56$p$481 &
		1,586122566110* &
		4961 &
		8161 &
		10\\
		\cellcolor{midgrau}
		&\cellcolor{hellgrau}
		15 &
		1,5625 &
		45 &
		75 &
		11\\
		\cellcolor{midgrau}
		&\cellcolor{hellgrau}
		13,225 &
		1,4894168402$p$7 &
		1679 &
		2929 &
		12\\
		\cellcolor{midgrau}
		&\cellcolor{hellgrau}
		12,25 &
		1,45001736$p$1 &
		161\textsuperscript{1} &
		289 &
		13\\
		\cellcolor{midgrau}
		&\cellcolor{hellgrau}
		11,75$p$5&
		1,43023882030* &
		1771 &
		3229 &
		14\\
		\cellcolor{midgrau}
		&\cellcolor{hellgrau}
		10,666$p$6 &
		1,3871604938* &
		[112]
		
		56&
		[212]
		
		106\textsuperscript{1} &
		15\\ \midrule
		\multicolumn{5}{m{14cm}}{* abbreviated number
			\ \ \ \textsuperscript{1} corrected number
			\ \ \ […] initial values prior shortening
			
			\textsuperscript{\$} for handling with decimal numbers we have to add 60 instead of 1. 
			
			The values in column I, were set by this way to a higher starting power so we have to multiply them by 60².
		
	$p$ within the decimal Values, indicates a periodic or repeating decimal number starting from that point.}
		\\\bottomrule
	\end{tabular}
\end{table}

\newpage
Likewise, not only the function value system is clearly demonstrable and by this way a veritable Confirmation of the correct translated Heading above, in addition a conscious conversion or transformation into measurable lengths can be proven with this deciphered calculation procedure.

That, this 3800 year old trigonometric function value system is unchanged from todays one and that our current values are exactly the same, can be shown in the following example as representative of all lines on Plimpton 322:

We look at line 3 of the cuneiform tablet: 
To work with the exact ratio of this table line, we simply take the decimal numerical value from column II and column III, this results in
\[
\frac{\textbf{4601}}{\textbf{6649}} = 0.69198375\footnote{abbreviated Number}
\]
The values to calculate can be read in the given order at the tablet. Their quotient or the ratio of Column II divided by Column III corresponds to the ``sine'' function value of the included angle $\alpha$. If you take the ``modern named angle'' from it, using a modern sine table or a computer, this corresponds to 43.787346\footnote{abbreviated Number} degrees. We can determine from this, the $\text{tan}\, \alpha = 0.958541p6$ and the value for the reciprocal function of cosine, the secant. 
\[
\text{sec}\, \alpha = 1 / \text{cos}\, \alpha = 1.385208p3.\footnote{$p$ means a periodic or repeating value from this point}
\]
So far we have done nothing else, than to use the given triangular proportions of the old Babylonian cuneiform tablet and to transform this triangle in our current trigonometric function value system.

Our value system relates to the base length in the unit circle of radius $1$. This was different in Babylonian times. We, therefore, simply multiply our current values by 60 and thus get the (Babylonian) Tangent and Secant values for line 3 of Plimpton 322 cuneiform tablet. 
\[
\text{bab-tan}\, \alpha = 57.5125 \text{ ;  and bab-sec}\, \alpha = 83.1125
\]

Column II and III can be interpreted as two sides of a right-angled triangle at first the opposite and thereafter the hypotenuse or diagonal. The decimal values of this line 3 are: 
\[
\textbf{4601};\ \textbf{6649}\ \text{;  and for the 3rd, not mentioned triangle side}\, \ b = \ \textbf{ 80} \times 60 = 4800.
\] 
If we now look at our bab-tan $\alpha$ = 57.5125 and bab-sec $\alpha$ = 83.1125 values as real fractions, we see them as:
\[
\text{bab-tan}\, \alpha =  \frac{\textbf{4601}}{\textbf{80}}, \ \text{bab-sec}\, \alpha  =  \frac{\textbf{6649}}{\textbf{80}}.
\]

By that, we get the triangle, documented at the cuneiform tablet Plimpton 322, line 3: 

\begin{center}
	Column II = \textbf{4601}\; Column III = \textbf{6649}\; ( The not mentioned 3rd Side: \textbf{80} $\times$ 60 = 4800),
\end{center}
direct without further calculations.
\\
\\
\\
From this point of view on, which is just as comprehensible for all other lines, it is clear, that we still use the same geometric construct for the trigonometric value system of today and in addition that the values from that time are just as exact as ours today.

\begin{landscape}
\begin{table}[p]
\centering
\caption{Calculation of Plimpton 322 with intermediate results with decimal numbers}
\begin{tiny}
\begin{tabular}{p{1.5cm}p{3cm}p{2.8cm}p{1cm}p{1.5cm}p{2.8cm}p{1.5cm}p{1.5cm}p{1.2cm}p{1.5cm}}
\toprule
\cellcolor{hellgrau}
	\textbf{Column 0} 

	(reconstr.)
	&
	\textbf{Column I}
	&
	\multicolumn{2}{l}{\cellcolor{hellblau}Reading and/or interim calculation} 
	&
	\textbf{Column II} 
	&
	\multicolumn{2}{l}{\cellcolor{hellblau}Reading and/or interim calculation}
	&
	\textbf{Column III}
	&
	\textbf{Column IV}
	&
	\cellcolor{midgrau}
	~
\\
\cellcolor{hellgrau}
	Note the doubled value, later add 1\textsuperscript{\$} and square this sum
	&
	The intermediate (takiltum) is the square of the diagonal, which minus 1 becomes the square of the width
	&
	\multicolumn{2}{p{4.3cm}}{\cellcolor{hellblau}Square root of Column I minus 1,
				 the Result is bab-tan $\alpha $ function value (to base 60)}
    &
	The Result for the width
	&
	\multicolumn{2}{p{4.8cm}}{\cellcolor{hellblau}Square root of Column I direct,
	the Result is bab-sec $\alpha $ function value (to base 60)}
	&
	The Result for the diagonal
	&
	It´s Place
	&
	\cellcolor{midgrau}
	Multiplier of the 3rd or base side (b) x 60
\\\midrule
\cellcolor{hellgrau}
	24,5
	&
	3600 x 1,983402p7 
	&
	\cellcolor{hellblau}59,5
		
	\cellcolor{hellblau}\textcolor{cyan}{\{59,5 x \textbf{2} = 119\}}
	&
	\cellcolor{hellblau}\centering  =$\frac{119} {\textbf{2}}$ 
	&
	119
	&
		\cellcolor{hellblau}
	84,5
		\cellcolor{hellblau}
	\textcolor{cyan}{\{84,5 x \textbf{2} = 169\}} &
		\cellcolor{hellblau}
	\centering  =$\frac{169} {\textbf{2}}$ 
	&
	169
	&
	1
	&
	\cellcolor{midgrau}\textbf{2} x 60 = 120	
\\
\cellcolor{hellgrau}
&
&\cellcolor{hellblau}
&\cellcolor{hellblau}
&
&\cellcolor{hellblau}
&\cellcolor{hellblau}
&
&
&\cellcolor{midgrau}
	\\
\cellcolor{hellgrau}
	23,76736p1 &
	3600 x 1,949158552088* &
		\cellcolor{hellblau}
	58,45486p1
		\cellcolor{hellblau}	
	\textcolor{cyan}{\{58,45486p1 x \textbf{288} = 16835\}} &
		\cellcolor{hellblau}
	\centering  =$\frac{16835}{\textbf{288}}$ &
	[16835] : 5
	
	=3367&
		\cellcolor{hellblau}
	83,76736p1
		\cellcolor{hellblau}
	\textcolor{cyan}{\{83,76736p1 x \textbf{288} = 24125\}} &
		\cellcolor{hellblau}
	\centering  =$\frac{24125}{\textbf{288}}$ &
	[24125] : 5
	
	=4825\textcolor{red}{\textsuperscript{1}} &
	2&
	\cellcolor{midgrau}
	[\textbf{288} x 60 ] : 5
	
	57,6 x60 = 3456
\\
\cellcolor{hellgrau}
&
&\cellcolor{hellblau}
&\cellcolor{hellblau}
&
&\cellcolor{hellblau}
&\cellcolor{hellblau}
&
&
&\cellcolor{midgrau}
	\\
	\cellcolor{hellgrau}
	23,1125&
	3600 x 1,918802126736p1&
		\cellcolor{hellblau}
	57,5125
	\cellcolor{hellblau}
	\textcolor{cyan}{\{57,5125 x \textbf{80} = 4601\}} &
		\cellcolor{hellblau}
	\centering  =$\frac{4601}{\textbf{80}}$ &
	4601&
		\cellcolor{hellblau}
	83,1125
	
		\cellcolor{hellblau}
	\textcolor{cyan}{\{83,1125 x \textbf{80} = 6649\}}&
		\cellcolor{hellblau}
	\centering  =$\frac{6649}{\textbf{80}}$ &
	6649 &
	3 &
	\cellcolor{midgrau}\textbf{80} x 60 = 4800
\\
\cellcolor{hellgrau}
&
&\cellcolor{hellblau}
&\cellcolor{hellblau}
&
&\cellcolor{hellblau}
&\cellcolor{hellblau}
&
&
&\cellcolor{midgrau}
	\\
	\cellcolor{hellgrau}
	22,40p4 &
	3600 x 1,886247906721* &

		\cellcolor{hellblau}
	56,48p4
			\cellcolor{hellblau}
	\textcolor{cyan}{\{56,48p4 x \textbf{225} 
		= 12709\}}&
		\cellcolor{hellblau}
	\centering  =$\frac{12709}{\textbf{225}}$ & 
	12709 &
		\cellcolor{hellblau}
	82,40p4
		\cellcolor{hellblau}
	\textcolor{cyan}{\{82,40p4 x \textbf{225} 
		= 18541\}} &
		\cellcolor{hellblau}
	\centering  =$\frac{18541}{\textbf{225}}$ &
	18541 &
	4&
	\cellcolor{midgrau}\textbf{225} x 60 = 13500
\\
\cellcolor{hellgrau}
&
&\cellcolor{hellblau}
&\cellcolor{hellblau}
&
&\cellcolor{hellblau}
&\cellcolor{hellblau}
&
&
&\cellcolor{midgrau}
	\\
	\cellcolor{hellgrau}
	20,8p3 &
	3600 x 1,815007716049* &
		\cellcolor{hellblau}
	54,1p6
		\cellcolor{hellblau}	
	\textcolor{cyan}{\{54,1p6 x \textbf{6} = 325\}}&
		\cellcolor{hellblau}
	\centering  =$\frac{325} {\textbf{6}}$ &
	[325] : 5
	
	=65&
	80,8p3
		\cellcolor{hellblau}
	\textcolor{cyan}{\{80,8p3 x \textbf{6} =485\}} &
		\cellcolor{hellblau}
	\centering  =$\frac{485} {\textbf{6}}$ &
	[485] : 5
	
	=97 &
	5 &
	\cellcolor{midgrau}[\textbf{6} x 60] : 5
	
	1,2 x 60 = 72
\\
\cellcolor{hellgrau}
&
&\cellcolor{hellblau}
&\cellcolor{hellblau}
&
&\cellcolor{hellblau}
&\cellcolor{hellblau}
&
&
&\cellcolor{midgrau}
	\\
\cellcolor{hellgrau}
	20,1p6 &
	3600 x 1,785192901234* &
	\cellcolor{hellblau}
	53,1p6
	
	\cellcolor{hellblau}
	\textcolor{cyan}{\{53,1p6 x \textbf{6} = 319\}} &
		\cellcolor{hellblau}
	\centering  =$\frac{319} {\textbf{6}}$ &
	319&
		\cellcolor{hellblau}
	80,1p6
	\cellcolor{hellblau}
	\textcolor{cyan}{\{80,1p6 x \textbf{6} =481\}} &
		\cellcolor{hellblau}
	\centering  =$\frac{481} {\textbf{6}}$ &
	481 &
	6&
	\cellcolor{midgrau}\textbf{6} x 60 = 360
\\
\cellcolor{hellgrau}
&
&\cellcolor{hellblau}
&\cellcolor{hellblau}
&
&\cellcolor{hellblau}
&\cellcolor{hellblau}
&
&
&\cellcolor{midgrau}
	\\
	\cellcolor{hellgrau}
	18,6p8 &
	3600 x 1,719983676268* &
\cellcolor{hellblau}
	50,9p1
	
\cellcolor{hellblau}
	\textcolor{cyan}{\{50,9p1 x \textbf{45} = 2291\}} &
		\cellcolor{hellblau}
	\centering  =$\frac{2291}{\textbf{45}}$ &
	2291&
		\cellcolor{hellblau}
	78,6p8
	
		\cellcolor{hellblau}
	\textcolor{cyan}{\{78,6p8 x \textbf{45} = 3541\}} &
	
		\cellcolor{hellblau}
	\centering  =$\frac{3541}{\textbf{45}}$ &
	3541 &
	7 &
	\cellcolor{midgrau}\textbf{45} x 60 = 2700
	\\
	\cellcolor{hellgrau}
	&
	&\cellcolor{hellblau}
	&\cellcolor{hellblau}
	&
	&\cellcolor{hellblau}
	&\cellcolor{hellblau}
	&
	&
	&\cellcolor{midgrau}
	\\
	\cellcolor{hellgrau}
	18,0625&
	3600 x 1,692709418402p7\textcolor{red}{\textsuperscript{1}} &
	\cellcolor{hellblau}
	49,9375
	
	\cellcolor{hellblau}	
	\textcolor{cyan}{\{49,9375 x \textbf{16} = 799\}} &
		\cellcolor{hellblau}
	\centering  =$\frac{799}{\textbf{16}}$ &
	799&
		\cellcolor{hellblau}
	78,0625
	
		\cellcolor{hellblau}
	\textcolor{cyan}{\{78,0625 x \textbf{16} = 1249\}} &
		\cellcolor{hellblau}
	\centering  =$\frac{1249}{\textbf{16}}$ &
	1249 &
	8 &
	\cellcolor{midgrau}\textbf{16} x 60 = 960
\\
\cellcolor{hellgrau}
&
&\cellcolor{hellblau}
&\cellcolor{hellblau}
&
&\cellcolor{hellblau}
&\cellcolor{hellblau}
&
&
&\cellcolor{midgrau}
	\\
	\cellcolor{hellgrau}
	16,9&
	3600 x 1,642669p4&
	48,1
		\cellcolor{hellblau}
		
	\textcolor{cyan}{\{48,1 x\textbf{10} =481\}} &
		\cellcolor{hellblau}
	\centering  =$\frac{481}{\textbf{10}}$ &
	481\textcolor{red}{\textsuperscript{1}} &
		\cellcolor{hellblau}
	76,9
	
		\cellcolor{hellblau}
	\textcolor{cyan}{\{76,9 x\textbf{10} = 769\}} &
		\cellcolor{hellblau}
	\centering  =$\frac{769}{\textbf{10}}$ &
	769&
	9 &
	\cellcolor{midgrau}\textbf{10} x 60 = 600
\\
\cellcolor{hellgrau}
&
&\cellcolor{hellblau}
&\cellcolor{hellblau}
&
&\cellcolor{hellblau}
&\cellcolor{hellblau}
&
&
&\cellcolor{midgrau}
	\\
	\cellcolor{hellgrau}
	15,56p481&
	3600 x 1,586122566110* &
	\cellcolor{hellblau}
	45,935p185
	\cellcolor{hellblau}
		
	\textcolor{cyan}{\{45,935p185 x \textbf{108} = 4961\}} &
	\cellcolor{hellblau}
	\centering  =$\frac{4961}{\textbf{108}}$ &
	4961&
		\cellcolor{hellblau}
	75,56p481
		\cellcolor{hellblau}
	
	\textcolor{cyan}{\{75,56p481 x \textbf{108} = 8161\}} &
		\cellcolor{hellblau}
	\centering  =$\frac{8161}{\textbf{108}}$ &
	8161&
	10 &
	\cellcolor{midgrau}\textbf{108} x 60 = 6480
\\
\cellcolor{hellgrau}
&
&\cellcolor{hellblau}
&\cellcolor{hellblau}
&
&\cellcolor{hellblau}
&\cellcolor{hellblau}
&
&
&\cellcolor{midgrau}
	\\
	\cellcolor{hellgrau}
	15 &
	3600 x 1,5625 &
		\cellcolor{hellblau}
	45
		\cellcolor{hellblau}
		
	\textcolor{cyan}{\{45 x \textbf{1} = 45\}} &
		\cellcolor{hellblau}
	\centering  =$\frac{45} {\textbf{1}}$ &
	45 &
	75
	\cellcolor{hellblau}
	
	\textcolor{cyan}{\{75 x \textbf{1} = 75\}} &
		\cellcolor{hellblau}
	\centering  =$\frac{75} {\textbf{1}}$ &
	75 &
	11 &
	\cellcolor{midgrau}\textbf{1} x 60 = 60
\\
\cellcolor{hellgrau}
&
&\cellcolor{hellblau}
&\cellcolor{hellblau}
&
&\cellcolor{hellblau}
&\cellcolor{hellblau}
&
&
&\cellcolor{midgrau}
	\\
\cellcolor{hellgrau}
	13,225 &
	3600 x 1,4894168402p7&
	\cellcolor{hellblau}
	41,975
	\cellcolor{hellblau}
		
	\textcolor{cyan}{\{41,975 x \textbf{40} = 1679\}} &
	\cellcolor{hellblau}
	\centering  =$\frac{1679}{\textbf{40}}$ &
	1679&
	\cellcolor{hellblau}
	73,225
	\cellcolor{hellblau}
	
	\textcolor{cyan}{\{73,225 x \textbf{40} = 2929\}} &
	\cellcolor{hellblau}
	\centering  =$\frac{2929}{\textbf{40}}$ &
	2929 &
	12 &
	\cellcolor{midgrau}\textbf{40} x 60 = 2400
\\
\cellcolor{hellgrau}
&
&\cellcolor{hellblau}
&\cellcolor{hellblau}
&
&\cellcolor{hellblau}
&\cellcolor{hellblau}
&
&
&\cellcolor{midgrau}
	\\
\cellcolor{hellgrau}
	12,25 &
	3600 x 1,45001736p1 &
	\cellcolor{hellblau}
	40,25
	\cellcolor{hellblau}
		
	\textcolor{cyan}{\{40,25 x \textbf{4} = 161\}} &
	\cellcolor{hellblau}
	\centering  =$\frac{161} {\textbf{4}}$ &
	\textcolor{black}{161}\textcolor{red}{\textsuperscript{1}} &
	\cellcolor{hellblau}
	72,25
	\cellcolor{hellblau}
	
	\textcolor{cyan}{\{72,25 x \textbf{4} = 289\}} &
	\cellcolor{hellblau}
	\centering  =$\frac{289} {\textbf{4}}$ &
	289&
	13 &
	\cellcolor{midgrau}\textbf{4} x 60 = 240
\\
\cellcolor{hellgrau}
&
&\cellcolor{hellblau}
&\cellcolor{hellblau}
&
&\cellcolor{hellblau}
&\cellcolor{hellblau}
&
&
&\cellcolor{midgrau}
	\\
\cellcolor{hellgrau}
	11,7p5 &
	3600 x 1,43023882030*&
	\cellcolor{hellblau}
	39,3p5
	\cellcolor{hellblau}
		
	\textcolor{cyan}{\{39,3p5 x \textbf{45} = 1771\}} &
	\cellcolor{hellblau}
	\centering  =$\frac{1771}{\textbf{45}}$ &
	1771 &
	\cellcolor{hellblau}
	71,7p5
	\cellcolor{hellblau}
	
	\textcolor{cyan}{\{71,7p5 x \textbf{45} = 3229\}} &
	\cellcolor{hellblau}
	\centering  =$\frac{3229}{\textbf{45}}$&
	3229&
	14 &
	\cellcolor{midgrau}\textbf{45} x 60 = 2700
	\\
	\cellcolor{hellgrau}
	&
	&\cellcolor{hellblau}
	&\cellcolor{hellblau}
	&
	&\cellcolor{hellblau}
	&\cellcolor{hellblau}
	&
	&
	&\cellcolor{midgrau}
	\\
\cellcolor{hellgrau}
	10,66p6&
	3600 x 1,3871604938*&
	\cellcolor{hellblau}
	37,33p3
	\cellcolor{hellblau}
		
	\textcolor{cyan}{\{37,33p3 x \textbf{3} = 112\}} &
	\cellcolor{hellblau}
	\centering =$\frac{112} {\textbf{3}}$ &
	[112] : 2
	
	=56 &
	\cellcolor{hellblau}
	70,66p6
	\cellcolor{hellblau}
	
	\textcolor{cyan}{\{70,66p6 x \textbf{3} = 212\}} &
	\cellcolor{hellblau}
	\centering  =$\frac{212} {\textbf{3}}$ &
	[212] : 2
	
	=106\textcolor{red}{\textsuperscript{1}}&
	15&
	\cellcolor{midgrau}[\textbf{3} x 60] : 2
	
	1,5 x 60 = 90
	\\
	\cellcolor{hellgrau}
	&
	&\cellcolor{hellblau}
	&\cellcolor{hellblau}
	&
	&\cellcolor{hellblau}
	&\cellcolor{hellblau}
	&
	&
	&\cellcolor{midgrau}
	\\\midrule
	\multicolumn{10}{p{20cm}}{
		\textbf{*} = abbreviated number;  
		\ \ \ \textbf{p} within the decimal Values, indicates a periodic or repeating decimal number starting from that point.
		
		\textcolor{red}{\textsuperscript{1}} = Corrected number; 
		
		[\textbf{…}] = initial values prior shortening;
		
		\textbf{lightblue enlightened cells} = Interim Calculation not on the cuneiform tablet, calculation just for explanation. 
		
	 \textsuperscript{\$}For working with decimal numbers, we have to add $60$ instead of $1$. So for Column I, the multiplier is $60^{2} = 3600$.  The function values at the interim calculation columns for tangent $\alpha $and secant $\alpha $ are the same as our modern 
		trigonometric values, only multiplied by $60$. 
		
		For line 15: the triple $28; 53; 45$ is also valid. \  \  \  \  \ ;\  \  \  A possible Column to the left of Column 0 is not shown at this table.
	}
	\\\bottomrule
\end{tabular}
\end{tiny}
\end{table}
\end{landscape}

\subsection{The calculation procedure step by step}

\begin{figure}[h]
	\centering
	\begin{minipage}{.5\textwidth}
		\centering
		\includegraphics[width=.92\textwidth]{img002.png}
		\caption{}
		\label{fig:abb1}
	\end{minipage}%
	\begin{minipage}{.5\textwidth}
		\centering
		\includegraphics[width=.9\textwidth]{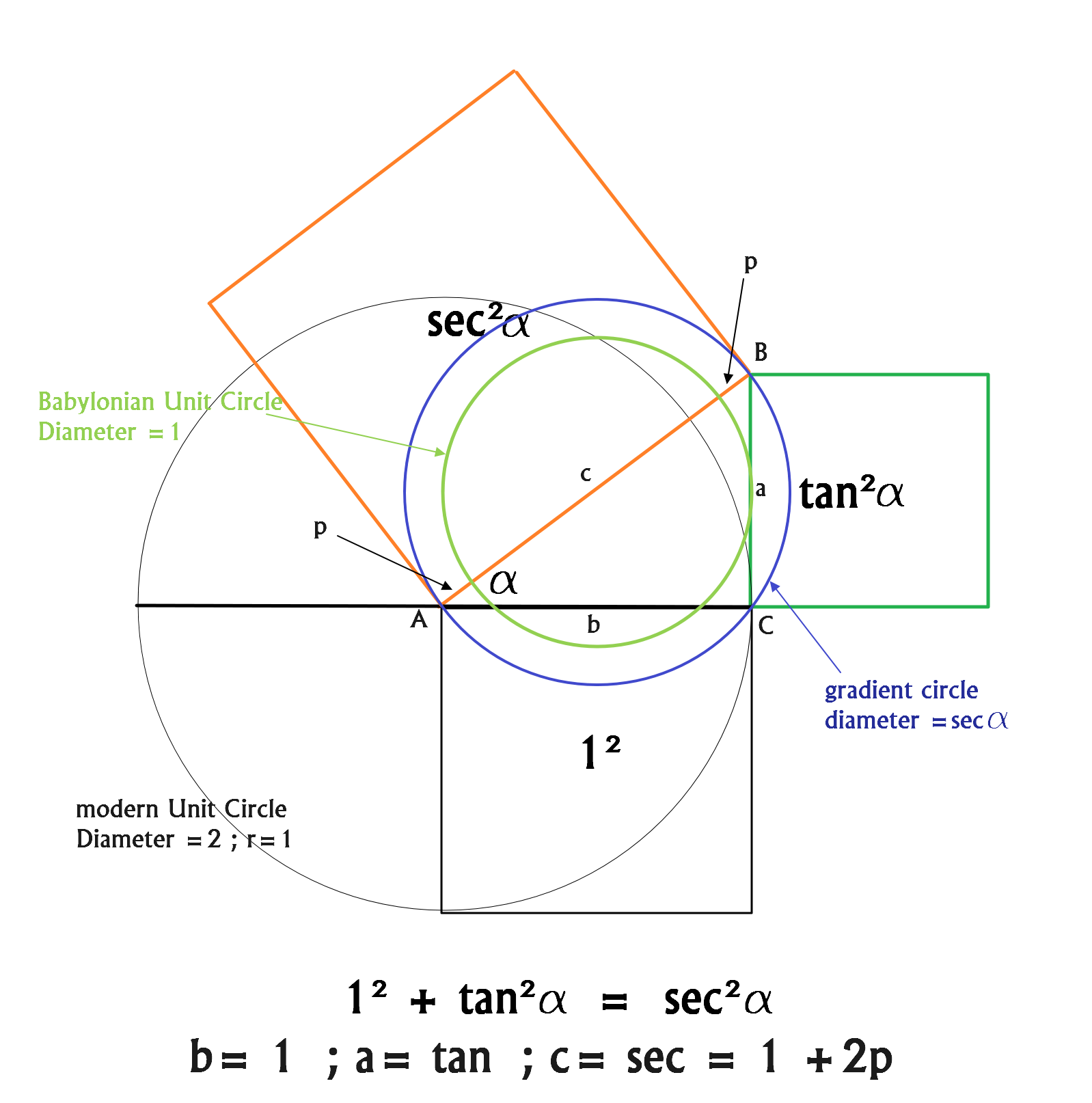}
		\caption{}
		\label{fig:abb2}
	\end{minipage}
\end{figure}
\label{section:stepbystep}
\begin{enumerate}
\item[a)]In column I is an intermediate result (Ta-ki-il-tum, Ta-ki-il-ti; pl.) given as a ratio or function value. The heading reads: \textit{The intermediate results are the squares of the diagonals (to the measured gradient ratio), this minus 1 becomes the square of the width (opening of the triangle).}
\[
	c^{2}- 1^{2} = a^{2} \text{ or rather } \text{sec}^{2}\alpha -1^{2} = \text{tan}^{2}\alpha.
\]
This heading corresponds to the so-called 2nd trigonometric Pythagoras (see figure \{\ref{fig:abb2}\}) and provides information on the further calculation sequence.

\item[b)]The result in column II is obtained from column I via an interim calculation:
Column I minus the value of sexagesimal 1 and take the square root of this. This result in the width (tangent) at the unit circle as a ratio or function value.  
\[
	\text{(base}=b=1\text{), therefore } \sqrt{\left(\dfrac{a}{b}\right)^{2}} = \dfrac{a}{b} = \dfrac{a}{1}
\]\[	
 \text{important is } b' = X\times 1 = X\times b\text{, so in form of a fraction } 
\] 
\[ b\text{ can be splitted into the fixed part $1$   
	and its stretching factor $Ext.$ or $X$.} 
\] 
\[ \text{by this, the function value or ratio, which is not already an integer, can be extended 
	by its fractional part to $1$:}
\]
\[
\dfrac{a}{b\times X}\times \dfrac{X}{1} = \dfrac{a'}{1} ; (X = Ext.)  \text{ , for example if $a$ is $4,25$ =}\dfrac{17}{4} \text{ this is: } \dfrac{(\dfrac{17}{4})}{1} \times \dfrac{4}{1} = \dfrac{17}{1\times 4} \times \dfrac{4}{1} = \dfrac{17}{1} = \dfrac {a'}{1}
\] 
\\                                                             
The ratio of the triangle sides $\dfrac{a}{b}$ or exact function value of tangent $\alpha$ was now extended to an integer Value (according to today's term) for further adoption and representability, by multiplication with the stretching factor $X$. In the example above it was $0,25 = \dfrac{1}{4} \text{ ; its reciprocal: }  \dfrac{4}{1} = X$. 

Here, the dimensionless function value, again became a measurable length $a'$  through the Extension. (for the ratios see also the headings at K.Vogel's Table on Plimpton 322\cite{vogel1959a},also\cite[98]{lehmann1992}).
 \,  \
 
 \[
 	a \times Ext. = a'\text{;}\text{ Extension = stretching factor $ X$ of the gradient triangle}
 \]
 \,  \

The above-mentioned square root values were read from tables with a very high degree of probability, so that this intermediate result did not have to be noted on the cuneiform tablet. At the same time and in typical smart Babylonian manner, any possible calculation- or spelling-error was filtered out and could not interfere the following steps. A confirmation for this one can see at line 8 \{\ref{section:errors}\}. 

The extension ($Ext. = X$) of the intermediate result and function value ``$a$'' ($sqrt$ of Column I minus $1$) is done in the decimal, directly by extension with the denominator of the fraction, if see the value as a real fraction. 
\[ 
   \text{(exemplary  at Line 1, for the width:   } a = 59,5 =\dfrac{119}{\textbf{2}}\text{ ;  } \dfrac{119}{2} \times  \textbf{2} = 119= a'
\]
\[    
   \text {  for its diagonal:   }  c = 84,5 =\dfrac{169}{\textbf{2}} ;   \text{ and }   \dfrac{169}{2} \times \textbf{2}  = 169 = c')
\]
\\
In the Babylonian sexagesimal system, this was done through an extension with the reciprocal, of a common factor in the value of $a$ and $c$. (also \cite{abdulaziz2010}) 

\[
(\text{exemplary   at Line 1:  } a= 59.30 ;  c= 1.24.30  ; \text{The common Factor is   }  30, \text{it's reciprocal  } \textbf{2}
\]
\[
\text{this result in  } 59.30\times \textbf{2} = 1.59 =a'  \text{   and  }  1.24.30 \times \textbf{2} = 2.49 = c' )
\]
\\
\item[c)]The now anticipated result in column III is created by taking the direct square root of the value of column I including the prefixed sexagesimal 1 (due to the squared value, this has the value 3600 = 60² as described in the table). This gives the ratio $\dfrac{c}{b}$ or the function value ``diagonal'' (the secant = the hypotenuse of the gradient triangle in the unit circle). 
Here, like in column II, the dimensionless value will be extended for further adoption and represent-ability with the same denominator as in column II, easily recognizable in the decimal fraction notation. In the Sexagesimal Notation is will be the same factor as in column II and its reciprocal to the extension. (See examples in the previous paragraph)

The reciprocal numerical value of the common factor contained in the ratio or function values is always the stretching factor \textit{X'} for all 3 sides of the respective triangle. The three sides of the triangle are for all gradient ratios or implicitly the angles $\alpha $ the following:
                   
     \[                      	  
             a' = X \times tan \, \alpha \,  \ \text{     ;  }  c' = X \times sec \, \alpha   \,  \ \text{     ;  } 
             b' = X \times 1  \text{  (  the Value for the Sexagesimal $1$ is $60$  )}
     \]

So one can see that the stretching $X$ is equal to the numerical value of the adjacent side $b'$ of the triangle in sexagesimal notation.
 
(Here the basic dimension 1 = 60 is not noted separately and we only see this multiplier or stretching factor.)

For this reason, the value of the 3rd side was not worth mentioning from the Babylonian point of view, since it was already included as a factor in the other two sides of the triangle.
For further explanation, refer to the tables attached here.

In addition to the last section c) of the procedure, the following can be specified: The Babylonian author knew the geometric construct on which his calculations were based very precisely and could read out this basic cathetus or 3rd page (the stretching factor) at any time if it was needed. In many lines of Plimpton 322, one can immediately subtract the portion of the value.
 
These portions of the diagonals or stretched secant values are known today as external secant or, in Babylonian, the arrow value, multiplied by the common factor, from the diagonal value in column III. 

If necessary, one received by this simple subtraction the stretching value of the base $b$, which always had to be an integral multiple of $60$.(see\footnote{Exemplary we can see this also at the YALE tablet ybc7289.    The Value for the diagonal is just an approximation,but it shows the principle: The Value is $42.25.35$ , the stretching factor is \textit{30}. The coefficient (\cite{fowler1998}), 
	
	\begin{center}
		or function value “diagonal = modern secant” is $1.24.51.10$.
	\end{center}
	
$42.25.35$ minus $30$ is the "stretched extern secant value": $12.25.35$. This divided by $30$ (the stretching factor) is: $24.51.10$. 
		
	This is the \textit{extern secant} portion, which was known in Babylonian time as $2\ \times $ \text{Arrow or: }$2p$. 
	So we return to Fig. 4:   \[ \textbf{1 + 2p = \text{Diagonal}  \text{   ;   } 1.0.0.0 + 24.51.10 =1.24.51.10}\] for our decimal system this means the sexagesimal \textit{24.51.10 + 1.0.0.0} is equal to the \textit{$\textbf{exsec} \text{ of } 45 \text{ Degree} $ + 1}, both is $\sqrt{2}$.} below)
\\ 
For a better overview, no subsequent shortening of the triangular values below base 1 = 60, was therefore initially considered. (See also: \{\ref{section:errors} \})

Furthermore, one would always have the option of checking via the
usual, but very rough, Babylonian approximation:
\[ 
\text{diagonal} \sim \text{base} + \dfrac{\text{width}^2 }{2 \times \text{base}}   ; \text{ \cite{lehmann1992} }
\]
\[
\text{ or via the correct calculation, by using :  }
	\text{   diagonal} = \sqrt{\text{base}^2 + \text{width}^2} \text{  (   Pythagorean theorem  ) }
\]
\[
\text{for calculation of the third side of the triangle.}
\]

Nevertheless, this base side was not necessary for the description and documentation of the gradient or further scale-independent use.
\\
Since the first analysis of the cuneiform tablet by O. Neugebauer and others (\cite{neugebauer1945};\cite{vogel1959a}), the reading from logic is conditionally a \textit{\underline{ratio value} at column I} and a \textit{triangle side for columns II and II}. 

The question of why exactly these triangle side lengths were noted, for the respective mentioned ratio value \\ in column I, could now be plausible explained for the first time.
\\

\item[d)]However, another question still arises as to why this geometric construct was chosen and what the initial numerical values of the respective line were?!
\\

For this we look at other Babylonian cuneiform artifacts in which parts of the triangle inscribed in it were calculated using a geometry in a circle and its diameter. It is noticeable that the Babylonians understood how to divide the diagonal = circle diameter, i.e. the hypotenuse of the triangle into a sum of the values 1 and 2P. The arrow p is the perpendicular line of the chord to the corresponding arc segment of the gradient circle. (See Fig. 1) In today's parlance, the arrow, also known as the rise of a bow or arch, is used to refer at the section of the bow delimited by this. Scriba \& Schreiber 2001 (\cite[18]{scriba2001}), rightly speaks of, as one can now be seen, "the beginning of the chord geometry", which was only made known much later, by Hipparchus. 

J. Lehmann 1994 (\cite[113]{lehmann1992}) at Exercise 36 also gives the same example in which both, the geometry in the circle and the principles later named after Pythagoras are incorporated.
 
With the same content also other cuneiform artifacts are dealing:  
\[2p + \text{base length of the triangle (b)} = \text{diagonal or diameter}\]
, this geometric construct is used several times and is even specifically linked to the descent, i.e. a negative slope (among others on BM 085194, CDLI No: P274707 exercise r i33 and r i39. (\cite{oracc2007}) 

This is not surprising either, because the arrow measure corresponds to one of the smallest function values in trigonometry, we know $"2p"$ in our current unit circle as $ Exsecans \alpha $.

In view of these examples, it is logical to assume this value before the given ratio of \textit{column I} and, as we now know, the function value of a unique gradient triangle, and it's circumcircle. 

It cannot be said yet, whether this was measured directly, or instead the appropriate arrow value for this gradient triangle was used via an auxiliary construction and measurement. 
However, it forms the logical transition to the value in \textit{column I}. 

However, one can imagine simple geometric constructs or devices with which the respective arrow measure could possibly also be determined directly.
\\

So you could set the following short section directly at the beginning of the now presented ``step by step'' description of the calculation sequence:

One determine the arrow value of the respective gradient triangle and write it down in column 0 (no longer exists today). It is plausible that the column heading looked like this: 
\begin{center}
	\textit{"write down the doubled value. Later add $1$ and create the square of this sum."}\footnote{ It needs to be explained in addition to the Babylonian calculus (see \{\ref{section:calculus}\}) that, the addition of $1$ leads to the diagonal or function value sec $\alpha$. 
		
	How this Number is to add for instance as, $1.0.0.0$ or just $1.0$ depends on the value of $2p$. What was surely very easy for the Babylonians, is sometimes hard to imagine for us today. Dependent from the inclination of the gradient triangle, the correct position at the sexagesimal number of the function value is to set, whether as the leading sexagesimal place of the sexagesimal number, or an addition to them. In our modern decimal system it is much easier, here is the Value for $2p$ equal to $exsec\,\alpha$.
	\[ \text{And so is:  }  exsec\,\alpha + 1 =sec\,\alpha, \text{thereby is $1$ anytime and easy to handle, the number : } 1,00 \] } 
\end{center}\ \
This then leads us directly to the values of \textit{column I}.
\\
\end{enumerate}

\begin{landscape}
\begin{table}[p]
\centering
\caption{Calculation of Plimpton 322 with sexagesimal numbers}
\begin{tiny}
	\begin{tabular}{m{2cm}m{3cm}m{3cm}m{1.5cm}m{3cm}m{2.5cm}m{2cm}m{1.2cm}m{1.5cm}}\toprule
		\cellcolor{hellgrau}{\bfseries Column 0 (reconstr.)} &
		{\bfseries Column I} &
		\cellcolor{hellblau}
		{\bfseries Reading and/or 
			
			interim calculation} &
		\cellcolor{midblau}
		\centering{\bfseries {\textless}{\textless}{\textless} {\textgreater}{\textgreater}{\textgreater}} &
		\cellcolor{hellblau}
		{\bfseries Reading and/or 
			
			interim calculation} &
		{\bfseries Column II} &
		{\bfseries Column III} &
		{\bfseries Column IV} &
		\cellcolor{midgrau}
		~
		\\
		\cellcolor{hellgrau}
		{}[Note the doubled value, later add 1 and square this sum] &
		The intermediate (takiltum) is the square of the diagonal, which minus 1 becomes the square of the width &
		\cellcolor{hellblau}
		Square root of Column I minus 1. 
		
		The interim result is bab-tan $\alpha $ function value multiplied with the reciprocal of the common factor &
		\cellcolor{midblau}
		The common
		 
		factor and 
		
		its \textbf{reciprocal} &
		\cellcolor{hellblau}
		Square root 
		of Column I direct. The interim result is bab-sec $\alpha $ function value multiplied with the 
		reciprocal of the common factor &
		The Result for the width &
		The Result for the diagonal &
		It´s Place &
		\cellcolor{midgrau}
		~
		
		Multiplier of the base side x 60\\\midrule
		\cellcolor{hellgrau}24.30  &
		1.59.00.15 &
		\cellcolor{hellblau}
		59.30
		\cellcolor{hellblau}
		
		\textcolor{cyan}{\{59.30 x \textbf{2} = 1.59\}} &
		\cellcolor{midblau}
		30
		\cellcolor{midblau}
		
		\textbf{2} &
		\cellcolor{hellblau}
		1.24.30
		\cellcolor{hellblau}
		
		\textcolor{cyan}{\{1.24.30 x \textbf{2} = 2.49\}}\ \  &
		1.59 &
		2.49 &
		1 &
		\cellcolor{midgrau}\textbf{2}\\
	\cellcolor{hellgrau}
	&
	&\cellcolor{hellblau}
	&\cellcolor{midblau}
	&\cellcolor{hellblau}
	&
	&
	&
	&\cellcolor{midgrau}
	\\
		\cellcolor{hellgrau}23.46.2.30 &
		1.56.56.58.14.50.06.15 &
		\cellcolor{hellblau}
		58.27.17.30
		\cellcolor{hellblau}
		
		\textcolor{cyan}{\{58.27.17.30 x \textbf{4.48} = 4.40.35\}} &
		\cellcolor{midblau}
		12.30
		\cellcolor{midblau}
		
		\textbf{4.48} &
		\cellcolor{hellblau}
		1.23.46.2.30
		\cellcolor{hellblau}
		
		\textcolor{cyan}{\{1.23.46.2.30 x \textbf{4.48} = 6.42.05\}} &
		[4.40.35] : 5
		
		56.07 &
		[6.42.05] : 5
		
		1.20.25 \textcolor{red}{\textsuperscript{1}}&
		2 &
		\cellcolor{midgrau}[\textbf{4.48}] : 5
		
		57.36\\
	\cellcolor{hellgrau}
	&
	&\cellcolor{hellblau}
	&\cellcolor{midblau}
	&\cellcolor{hellblau}
	&
	&
	&
	&\cellcolor{midgrau}
	\\
		\cellcolor{hellgrau}23.6.45 &
		1.55.07.41.15.33.45 &
		\cellcolor{hellblau}
		57.30.45
		\cellcolor{hellblau}
		
		\textcolor{cyan}{\{57.30.45 x \textbf{1.20} = 1.16.41\}} &
		\cellcolor{midblau}
		45
		\cellcolor{midblau}
		
		\textbf{1.20} &
		\cellcolor{hellblau}
		1.23.06.45
		\cellcolor{hellblau}
		
		\textcolor{cyan}{\{1.23.06.45 x \textbf{1.20} = 1.50.49\}} &
		1.16.41 &
		1.50.49 &
		3 &
		\cellcolor{midgrau}\textbf{1.20}\\
	\cellcolor{hellgrau}
	&
	&\cellcolor{hellblau}
	&\cellcolor{midblau}
	&\cellcolor{hellblau}
	&
	&
	&
	&\cellcolor{midgrau}
	\\
		\cellcolor{hellgrau}22.24.16 &
		1.53.10.29.32.52.16 &
		\cellcolor{hellblau}
		56.29.04
		\cellcolor{hellblau}
		
		\textcolor{cyan}{\{56.29.04 x \textbf{3.45} = 3.31.49\}} &
		\cellcolor{midblau}
		16
		\cellcolor{midblau}
		
		\textbf{3.45} &
		\cellcolor{hellblau}
		1.22.24.16
		\cellcolor{hellblau}
		
		\textcolor{cyan}{\{1.22.24.16 x \textbf{3.45}= 5.9.1\}} &
		3.31.49 &
		5.09.01 &
		4 &
		\cellcolor{midgrau}\textbf{3.45}\\
		\cellcolor{hellgrau}
	&
	&\cellcolor{hellblau}
	&\cellcolor{midblau}
	&\cellcolor{hellblau}
	&
	&
	&
	&\cellcolor{midgrau}
	\\
		\cellcolor{hellgrau}20.50 &
		1.48.54.01.40 &
		\cellcolor{hellblau}
		54.10
		\cellcolor{hellblau}
		
		\textcolor{cyan}{\{54.10 x \textbf{6} = 5.25\}} &
		\cellcolor{midblau}
		10
		\cellcolor{midblau}
		
		\textbf{6} &
		\cellcolor{hellblau}
		1.20.50
		\cellcolor{hellblau}
		
		\textcolor{cyan}{\{1.20.50 x \textbf{6} = 8.05\}} &
		[5.25] : 5
		
		1.05 &
		[8.05] : 5
		
		1.37 &
		5 &
		\cellcolor{midgrau}[\textbf{6} ] : 5
		
		1.12\\
	\cellcolor{hellgrau}
	&
	&\cellcolor{hellblau}
	&\cellcolor{midblau}
	&\cellcolor{hellblau}
	&
	&
	&
	&\cellcolor{midgrau}
	\\
		\cellcolor{hellgrau}20.10 &
		1.47.06.41.40 &
		\cellcolor{hellblau}
		53.10
		\cellcolor{hellblau}
		
		\textcolor{cyan}{\{53.10 x \textbf{6} = 5.19\}} &
		\cellcolor{midblau}
		10
		\cellcolor{midblau}
		
		\textbf{6} &
		\cellcolor{hellblau}
		1.20.10
		\cellcolor{hellblau}
		
		\textcolor{cyan}{\{1.20.10 x \textbf{6} = 8.01\}} &
		5.19 &
		8.01 &
		6 &
		\cellcolor{midgrau}\textbf{6}\\
	\cellcolor{hellgrau}
	&
	&\cellcolor{hellblau}
	&\cellcolor{midblau}
	&\cellcolor{hellblau}
	&
	&
	&
	&\cellcolor{midgrau}
	\\
		\cellcolor{hellgrau}18.41.20 &
		1.43.11.56.28.26.40 &
		\cellcolor{hellblau}
		50.54.40
		\cellcolor{hellblau}
		\textcolor{cyan}{\{50.54.40 x \textbf{45} 
			= 38.11\}} &
		\cellcolor{midblau}
		1.20 
		\cellcolor{midblau}
		
		\textbf{45} &
		\cellcolor{hellblau}
		1.18.41.20
		\cellcolor{hellblau}
		\textcolor{cyan}{\{1.18.41.20 x \textbf{45} 
			= 59.01\}} &
		38.11 &
		59.01 &
		7 &
		\cellcolor{midgrau}\textbf{45}\\
		\cellcolor{hellgrau}
		&
		&\cellcolor{hellblau}
		&\cellcolor{midblau}
		&\cellcolor{hellblau}
		&
		&
		&
		&\cellcolor{midgrau}
		\\
		\cellcolor{hellgrau}18.3.45 &
		1.41.33.45.14.03.45\textcolor{red}{\textsuperscript{1}} &
		\cellcolor{hellblau}
		49.56.15
		\cellcolor{hellblau}
		
		\textcolor{cyan}{\{49.56.15 x \textbf{16} = 13.19\}} &
		\cellcolor{midblau}
		3.45
		\cellcolor{midblau}
		
		\textbf{16} &
		\cellcolor{hellblau}
		1.18.3.45
		\cellcolor{hellblau}
		
		\textcolor{cyan}{\{1.18.3.45 x \textbf{16} = 20.49\}} &
		13.19 &
		20.49 &
		8 &
		\cellcolor{midgrau}\textbf{16}\\
		\cellcolor{hellgrau}
		&
		&\cellcolor{hellblau}
		&\cellcolor{midblau}
		&\cellcolor{hellblau}
		&
		&
		&
		&\cellcolor{midgrau}
		\\
		\cellcolor{hellgrau}16.54 &
		1.38.33.36.36 &
		\cellcolor{hellblau}
		48.06
		\cellcolor{hellblau}
		
		\textcolor{cyan}{\{48.06 x\textbf{10} =8.01\}} &
		\cellcolor{midblau}
		6
		\cellcolor{midblau}
		
		\textbf{10} &
		\cellcolor{hellblau}
		1.16.54
		\cellcolor{hellblau}
		
		\textcolor{cyan}{\{1.16.54 x \textbf{10} = 12.49\}} &
		8.01\textcolor{red}{\textsuperscript{1}} &
		12.49 &
		9 &
		\cellcolor{midgrau}\textbf{10}\\
		\cellcolor{hellgrau}
		&
		&\cellcolor{hellblau}
		&\cellcolor{midblau}
		&\cellcolor{hellblau}
		&
		&
		&
		&\cellcolor{midgrau}
		\\
		\cellcolor{hellgrau}15.33.53.20 &
		1.35.10.02.28.27.24.26.40 &
		\cellcolor{hellblau}
		45.56.6.40
		\cellcolor{hellblau}
		
		\textcolor{cyan}{\{45.56.6.40 x \textbf{1.48} = 1.22.41\}} &
		\cellcolor{midblau}
		33.20
		\cellcolor{midblau}
		
		\textbf{1.48} &
		\cellcolor{hellblau}
		1.15.33.53.20
		\cellcolor{hellblau}
		
		\textcolor{cyan}{\{1.15.33.53.20 x \textbf{1.48} = 2.16.1\}} &
		1.22.41 &
		2.16. 1 &
		10 &
		\cellcolor{midgrau}\textbf{1.48}\\
		\cellcolor{hellgrau}
		&
		&\cellcolor{hellblau}
		&\cellcolor{midblau}
		&\cellcolor{hellblau}
		&
		&
		&
		&\cellcolor{midgrau}
		\\
		\cellcolor{hellgrau}15 &
		1.33.45 &
		\cellcolor{hellblau}
		45
		\cellcolor{hellblau}
		
		\textcolor{cyan}{\{45 x \textbf{1} = 45\}} &
		\cellcolor{midblau}
		1 
		\cellcolor{midblau}
		
		\textbf{1} &
		\cellcolor{hellblau}
		1.15
		\cellcolor{hellblau}
		
		\textcolor{cyan}{\{1.15 x \textbf{1} = 1.15\}} &
		45 &
		1.15 &
		11 &
		\cellcolor{midgrau}\textbf{1}\\
		\cellcolor{hellgrau}
		&
		&\cellcolor{hellblau}
		&\cellcolor{midblau}
		&\cellcolor{hellblau}
		&
		&
		&
		&\cellcolor{midgrau}
		\\
		\cellcolor{hellgrau}13.13.30 &
		1.29.21.54.02.15 &
		\cellcolor{hellblau}
		41.58.30
		\cellcolor{hellblau}
		
		\textcolor{cyan}{\{41.58.30 x \textbf{40} = 27.59\}} &
		\cellcolor{midblau}
		1.30
		\cellcolor{midblau}
		
		\textbf{40} &
		\cellcolor{hellblau}
		1.13.13.30
		\cellcolor{hellblau}
		
		\textcolor{cyan}{\{1.13.13.30 x \textbf{40} = 48.49\}} &
		27.59 &
		48.49 &
		12 &
		\cellcolor{midgrau}\textbf{40}\\
		\cellcolor{hellgrau}
		&
		&\cellcolor{hellblau}
		&\cellcolor{midblau}
		&\cellcolor{hellblau}
		&
		&
		&
		&\cellcolor{midgrau}
		\\
		\cellcolor{hellgrau}12.15 &
		1.27.00.03.45 &
		\cellcolor{hellblau}
		40.15
		\cellcolor{hellblau}
		
		\textcolor{cyan}{\{40.15 x \textbf{4} = 2.41\}} &
		\cellcolor{midblau}
		15
		\cellcolor{midblau}
		
		\textbf{4} &
		\cellcolor{hellblau}
		1.12.15
		\cellcolor{hellblau}
		
		\textcolor{cyan}{\{1.12.15 x \textbf{4} = 4.49\}} &
		2.41 \textcolor{red}{\textsuperscript{1}}&
		4.49 &
		13 &
		\cellcolor{midgrau}\textbf{4}\\
		\cellcolor{hellgrau}
		&
		&\cellcolor{hellblau}
		&\cellcolor{midblau}
		&\cellcolor{hellblau}
		&
		&
		&
		&\cellcolor{midgrau}
		\\
		\cellcolor{hellgrau}11.45.20 &
		1.25.48.51.35.06.40 &
		\cellcolor{hellblau}
		39.21.20
		\cellcolor{hellblau}
		
		\textcolor{cyan}{\{39.21.20 x \textbf{45} = 29.31\}} &
		\cellcolor{midblau}
		1.20
		\cellcolor{midblau}
		
		\textbf{45} &
		\cellcolor{hellblau}
		1.11.45.20
		\cellcolor{hellblau}
		
		\textcolor{cyan}{\{1.11.45.20 x \textbf{45} = 53.49\}} &
		29.31 &
		53.49 &
		14 &
		\cellcolor{midgrau}\textbf{45}\\
		\cellcolor{hellgrau}
		&
		&\cellcolor{hellblau}
		&\cellcolor{midblau}
		&\cellcolor{hellblau}
		&
		&
		&
		&\cellcolor{midgrau}
		\\
		\cellcolor{hellgrau}10.40 &
		1.23.13.46.40 &
		\cellcolor{hellblau}
		37.20
		\cellcolor{hellblau}
		
		\textcolor{cyan}{\{37.20 x \textbf{3} = 1.52\}} &
		\cellcolor{midblau}
		20
		\cellcolor{midblau}
		
		\textbf{3} &
		\cellcolor{hellblau}
		1.10.40
		\cellcolor{hellblau}
		
		\textcolor{cyan}{\{1.10.40 x \textbf{3} = 3.32\}} &
		[1.52] : 2
		
		56 &
		[3.32] : 2
		
		1.46 \textcolor{red}{\textsuperscript{1}}&
		15 &
		\cellcolor{midgrau}[\textbf{3} ] : 2
		
		1.30\\
		\cellcolor{hellgrau}
		&
		&\cellcolor{hellblau}
		&\cellcolor{midblau}
		&\cellcolor{hellblau}
		&
		&
		&
		&\cellcolor{midgrau}
		\\\bottomrule
		\multicolumn{9}{p{20cm}}{
			\textit{Lightblue enlightened cells} = Interim Calculation not on the cuneiform tablet, calculation just for explanation \ \ ; \ \ \textit{white back-grounded cells} = visible Content of the Plimpton 322 cuneiform Tablet.
			
			\textcolor{red}{\textsuperscript{1}} = Corrected number ;   
			
			[…] = initial values prior shortening;
			
			\textit{Lightgrey enlightened cells} = Reconstructed Column 0 (left) , and the 3rd Triangle Side (right column just for explanation)

			The function values for tangent $\alpha $and secant $\alpha $ are the same as our values, multiplied by 60.

			For line 15: the triple 28; 53; 45 is also valid. /  A possible Column to the left of Column 0 is not shown at this table.
		}
	\\\bottomrule
	\end{tabular}
\end{tiny}
\end{table}
\end{landscape}

\subsection{The resulting boundary conditions of the calculation.}

The following boundary conditions result from the reconstructed calculation procedure:

I. The definition range of the calculation sequence is valid according to the given heading in
column I for right-angled triangles with a slope angle $\alpha $ between 0° and $\sim$ 82.55° \cite{Kleb2021a}.
Outside this range, the column heading is not consistent anymore, because $1$ has to be subtracted from the second sexagesimal place in order to get the "square of the width". 

II. For the correct gradient triangle and its proportions, column III must also be extended with the same factor (reciprocal of the common factor = stretching factor $X$) from column II (see \{\ref{section:errors}\}).

III. The results in columns II and III can be shortened. By shortening them afterwards, it is important to ensure that in this case for both values of column II and III, it must be shortened by the same number (reciprocal of the common factor).

A shortening is only permissible by 2 or 4 if both side values of the triangle have calculated end
places with a straight sexagesimal value.

A shortening by 3 is only permitted if the last sexagesimal place of both sides ending on 3 or
7, or in few cases one side is already an integer (also alternating).

A shortening by 5 or possibly its multiple values, is only permitted if the last sexagesimal
place of both values ends at 5 or full tens (see \{\ref{section:errors}\}).

IV. If the intermediate values for the columns II and III (which are calculated by taking the square-root) will not have a sexagesimal place after the 1st or 2nd place, the common extension factor is $1$. This is already visible at our decimal notation for these values, they are already integer values.

V. Since the sexagesimal notation knows no real fractions in today's sense, we instead multiply with the reciprocal of the common factor of the values instead of the denominator. 
This practice of multiplying with the reciprocal instead of dividing is not only common, there are also countless other examples for this procedure\cite{neugebauer1969}\cite{neugebauer1935}\cite{neugebauer1937}\cite{oracc2007}.
\\

\subsection{A Brief summary of the calculation procedure}

In addition, to the above-mentioned boundary condition III., the following statement can be made: The shortening of lines 2, 5 and 15 is highly plausible a separate step, which were proceeded after the extension with their unique reciprocal factor. 
Within line 15, for example, it can be clearly seen that the actual common factor of the function values $37.20$ and $1.10.40$ is $20$, this reciprocally is $3$. ($3 \times .20$ = $1$ ; and $3 \times .40$ = $2$  , both are integer)

The complete resulting values are $1.52$ and $3.32$, but after this calculation they were incorrectly shortened. Once divided by $2$ and once by $4$. This results in the incorrect assembled Values of $56$ for Column II and $53$ for Column III.

On the other hand it is clearly to understand that, the choice of two separate and initial chosen, but then `` not common'' factors and their reciprocals , is to exclude according to the normal logic. 
The therefore necessary extension factors for the valid function values $37.20$ and $1.10.40$ would be in this case for: $\text{Column II } 1.30 \text{ reciprocal } 40 \text{; for Column III otherwise } 45 \text{ reciprocal } 1.20$, but these selections is not plausible in any way.

For the first time, the complete calculation procedure of Plimpton 322 could be presented consistently. Accordingly, we are dealing with the diagonal calculation used by the Babylonians with a dimensionless
ratio or function value system, which can only be used in reality by multiplication, with a real existing measure of length. There are already further references to this in the literature for other cuneiform tablets that have been found, which can now be interpreted even better or could possible embedded in the overall context of these function value system (YBC7289; SI.427 etc.).

Even if the Babylonian Creators have used this trigonometric function system in a slightly different way and has not used the sine function value for that, we have even today multiplying a real length by the function value of sine, cosine, tangent, cotangent and secant too, to get the respective  real side length of the triangle. 
In order to describe an ascending angle today via a ratio, coefficient or function value which is all the same in meaning, we mostly use the equally dimensionless sine function. Therefore, it is only a marginal note, that Plimpton 322, presents us in column II ( the opposite side) and after that, column III ( the diagonal or hypotenuse), the necessary values to calculate this function value too. These are exactly the triangle sides which we need today for our sine function values result. (sin$\alpha $ = opposite / hypotenuse).

Even if the Babylonians don't use the sine, there is a common ground for all of this. Moreover, neither the Babylonian's, nor we today, need for this result the base side of the triangle. (3rd triangle side).

Despite, of the self-given limitations in the count of possible steps or valid gradient triangles, due to the special selected and computable Values to a minimum of 150 for a quarter circle,  (we can assume a regular N-corner with approx. 600 sides for the full circle) 
\begin{center}
	``the Babylonian diagonal calculation'',
\end{center}
is clearly to be equated with our trigonometric angle functions today, or correct spelled, to be regarded as its real predecessor.
\\
\\
\subsection{Analysis of errors on Plimpton 322}
\label{section:errors}

Over the years, many scientists have devoted themselves to the various errors on the
cuneiform tablet. Some of them could already be explained plausibly, for other errors this did not succeed. Scriba and Schreiber (\cite{scriba2001}p21) summarized this for line 2 of the cuneiform tablet as follows: (translation from German) \textit{An explanation of this error could perhaps help to uncover the construction principle.}
\\

\paragraph{Line 2}
The noted value in column III is: $3.12.1$, the corrected should be initially $6.42.05$. after shortening by $5$ it is $1.20.25$. This error is also the most complex, but it can be plausible deciphered here, that according to the Babylonian calculation, the original values for line 2 calculated by the author of the tablet were: in column II the width, in the decimal $16835$ or sexagesimal $4.40.35$, and for column III the diagonals, of the same line $24125$ or $6.42.05$.

Only now, and as a really unnecessary and additional last step, which was only possible and not prohibited in a few lines, were these values shortened again. However, certain rules had to be adhered here, if both values (width and diagonal) ended on $5$, this could only be shortened by dividing with $5$ or, if necessary, a multiple of $5$, according to the conditions.

Unfortunately, we have not the short intermediate calculation that was made between columns I and II because it was made as a mental calculation, or by means of another temporarily used writing board. Due to the given last sexagesimal place of 60th at the Values, there was only the possibility of shortening it by factor $5$, accordingly the value of $a$ or the width in column II was correct divided by $5$.

It must be said again as an important addition that, it can also be demonstrated in other lines that the lines was possibly calculated column per column with a certain time offset and independently of each other, or more plausible, line by line was calculated, but also with a lag of time between. However, this posed the risk that the numerical value directly next to it within the other column of the line, could not serve as a reference or control point, because it was not, or not yet, available.

For the diagonal value of column III, line 2, can be assumed that, the correct, unabridged length value $6.42.5$ was mentally held, or was written on an extern interim table of the interim calculation. However, perhaps due to poor lighting, or unclear notation the real final Number $5$, was interpreted as a $2$, because of the similarity of the cuneiform sign, who was also drawn with 2 vertical lines. Both containing 2 verticals and sometimes we can find cuneiform texts with a $4$, or a $5$, or a $6$ very narrow and with low distance vertical drawn, so that they look like $1$, $2$ or $3$. 

In this consequence the following misread numerical value to shorten, $6.42."2"$ was used, instead of the correct one above.

As a result of this misreading, a known rule was used for the value: With an even final number could be shortened, by division with $2$, or possibly also by $4$. 

All in all and in any way, we are dealing with an error of at least 2 levels here, since the existence of a second hook of a 10th at this place( Line 2 Column III), is not certain, are the next two explanations the most plausible ones.

Because a shortening by a division by $2$ was surely carried out safely in the head, there are only these logical explanations for the misread Value of ($6.42."2"$).

Either he mentally twisted the already distorted number $6.42."2"$ to $6."24"."2"$ and shortened this by $2$ to $3.12.1$. Such a mental twisted number in the head, knows even everyone who is very familiar and is working daily with mathematics.

Or secondly, the scribe shortened $6.42."2"$ correctly and received $3.21.1$, which he then noted incorrectly in the second step.

Alternatively, a double incorrect reading leads to $6."44"."2"$ would also be conceivable. $6."44"."2"$, was correctly divided by 2. This would result in the sexagesimal value of 3.22.1. Since line 2 of the cuneiform table is the narrowest line of the entire table and at least on the currently available images a second indistinct 10th hook in front of the other one cannot be ruled out to 100\%, this possible explanation should also be mentioned and discussed.

A final explanation, just for completion at least to be addressed here, can be ruled out. If the author had halved the correct value $6.42.5$ in spite of the fact that, it was not allowed to shorten by division by $2$, the result would have been $3.21.2.30$. In order to get to the noted value $3.12.1$, but more mistakes would have to be made for this, than in the two ( alternatively three) explanations before.

At this point we can conclude, the mentioned solutions can show us how the error was created and even this is very important, what was the thinking of the Babylonian author during the calculation steps. But, we have the real fact, that the initial value of Line 2 Column III was a probably misreading of the correct $6.42.5$.

To summarize correct values of the gradient triangle of line 2, they are: Column II $4.40.35 \, (16835)$; Column III $6.42.5 \, (24125)$; and the unmentioned 3rd Side $4.48 \,(288 \times 60 = 17280)$, only correct shortened by $5$, this results in $56.07 \, (3367)$; $1.20.25 \, (4825)$; and for the 3rd Side $57.36 \,(3456)$.
\\

I do not want to go into the alleged error in column I here, since it is “only” a too narrow place between the 6th and 7th sexagesimal place. The Babylonian author would have noticed an incorrect merging of $50$ and $6$ in the immediately following calculation step of extracting the roots or reading from a table of roots.

\paragraph{Line 5}
The initial values prior shortening by $5$, were for column II $5.25$ and column III $8.05$.
There are no errors in this line, but it should be noted that, the values
on the cuneiform table at line 5, are already the correct shortened Values by dividing with $5$.

\paragraph{Line 8}
The noted value in column I is $1.41.33.59.03.45$, the correct value is  $1.41.33.45.14.03.45$.
Line 8 contains an error in column I, which, however, did not affect the correctness of the other values; this confirms the automatic correction, which has already been mentioned several times, through the use of external square root value tables. Whether the error was caused by a mistake when manually squaring column 0, to column I or by copying incorrectly from a table cannot be clearly determined, but the first variant appears to be the more plausible.

\paragraph{Line 9}
The noted value in column II is $9.01$ instead of the correct $8.01$.
This is a clear spelling mistake, once set the author did not want to revise the wrong one, especially since, as already described, the several mathematical controls could have been used in case of further usage.

\paragraph{Line 13}
The noted value in column II is $7.12.1$. The correct value is $2.41$. 
Since the table was evidently created at a time when the cuneiform character for double zero, or as a placeholder was not used, the representation is correct with a large gap between the 2nd and 4th sexagesimal places in column I.

In spite of this, the error that follows in the II column is the most interesting. This confirms the calculation procedure presented here, again. 
The Babylonian took the initial value from column I for column II (minus $1.0.0.0$), or looked for it correctly in the table of roots. It can only be assumed that, through carelessness or distraction he got into the wrong column of the adjacent numerical values and probably immediately recognized a common factor of $3.45$ who was very familiar to him in the context of these calculations. Therefore, the following extension step was carried out from the actual initial value and not its root value! The author has multiplied the complete value of column I (minus $1.0.0.0$ for the value of column II) by the reciprocal of factor $3.45$. Its reciprocal is $16$. This results in $27.00.3.45 \times 16 = 7.12.1$.

t is also clear that the Babylonian author has plausibly not reached the actually correct function Value and then the extended value of the side length. On the one hand, because that, he would immediately notice that, the correct and already very compact result, the function value of bab-tan$\alpha $ = $40.15$ leads to its extension with the reciprocal of the common factor $.15$, this is $4$. 
On the other hand, no new squaring was planned during the calculation procedure, therefore it is not plausible to assume that, the Babylonian author has unnecessarily squared the correct $2.41$. 
In any way, the correct result of Line 13 is a triple of: $2.41; 4.49; (4.00)$.

\paragraph{Line 15}
The noted value in column II is $56$ and the value in column III is $53$. The correct values are $56$ and $1.46$ or alternatively $28$ and $53$.
The error in line 15 underlines once more the column by column or at minimum time lagged line by line processing, compared to the ongoing line by line step by step processing, that is actually normal for us. From the initial values in column I, the square root is the ratio $\frac{a}{b}$ = bab-tan$\,\alpha $ = $37.20$ for column II arose in an interim calculation or reading from the table of roots and the ratio $\frac{c}{b}$= bab-sec$\,\alpha $ = $1.10.40$ for column III.

The common factor $20$ and its reciprocal $3$, as the expansion multiplier $X$, results in a triple of $1.52; 3.32; (3.00)$. Nevertheless, in the following, the result of the columns was again shortened, so that column II corresponds to $1.52$ by $2$ = $56$. The correct result of the diagonal: $3.32$, however, was shortened by $4$, resulting in $53$.

At this point it can only be assumed that the original intention was not to record triples with a base value (adjacent) below $1 = 60$. From this sight, the triple $56; 1.46; (1.30)$ seems originally intended , even if the triple of $28; 53; (45)$ would also be computationally correct (see also Neugebauer and Sachs \cite[p. 38]{neugebauer1945}).

\paragraph{Summary}
In summary of the previous error analysis, the following can be stated. Although all the made errors, do not lead to the failure of an entire value line, whether at the time of creation or even today, it is obvious that many of the incorrect values are directly or indirectly related to the time-shifted processing of the tablet columns II and III.

There are only two possible scenarios: The Babylonian author was primarily concerned with fixing the gradient ratio by means of columns 0 and I from lines 1 to 15. Only then were the two columns II and III calculated independently of each other, i.e. column by column.

However, it seems more plausible and this supports the calculation procedure presented here, that the Babylonian mathematicians, who are known to be very efficient, has been used the already calculated result of the ratio or the bab-sec $\alpha $ at the transition from column 0 to column I, for the column III. 
This value was extended and noted, a few interim calculations ahead those of column II (because this Value could have been only calculated From the squared Secant Value at Column I, and after the subtraction of $1$ and the following taking of the square root). 
So, this would be a staggered / time lagged,  but much more plausible line-by-line processing of the values. But by this way the direct comparative reference to the value of the other column and its common factor, was missing.
\\
\\
\\

\subsection{Further conclusions and answers}

The calculation procedure gives rise to further questions, which are to be answered below.

\, \

\subsubsection{Is it an anachronism to suggest that, the Babylonians has used geometrics (inscribed triangles) within a circle and additional function values instead of real length values, at that time?}

As often mentioned within this work, there are many indications, that this was a normal doing and thinking during this time. So Scriba, Schreiber, Lehmann, and others has presented tasks and solutions from cuneiform tablets from this time. They show us the partition of the Diagonal, which is the hypotenuse of the triangle in 1+2p also(\cite[113]{lehmann1992}). This Partition is only be valid and makes sense with the circumcircle around the Triangle. 
The Oracc Database offers several Cuneiform artifacts who are dealing with Triangles, his doubled geometric Variant the rectangle and its connected circumcircle. (\cite{oracc2007} for instance BM085194)
Also in Egypt they have used the here discussed function values but in a kind of other way and with a separate mentioned Name the "Seked" \cite{bartclair2012}.
 
In conclusion, it is not an anachronism quite the contrary, it is strongly a part of the old Babylonian Mathematics and one of the normal ways to find a solution for some tasks.

\, \

\subsubsection{If it is a complete function value system, are there also other triples according to the calculation procedure used on Plimpton 322?}

Yes, far beyond the range shown on Plimpton 322, or at the discussed 38 lines in previous scientific literature, there are demonstrably more than 230 additional lines with valid gradient triangles. All of them are within the given definition range by the Heading of Column I \cite{kleb2020}\cite{Kleb2021a} and so between 0 and 82°. Due to the number of sexagesimal places, one can assume that only a subset of it was used, but that there were at least 150 proportions or function values in the mentioned range which could be used.

It can also be probably understood, that approximations were also included for important ratios, such as 1: 1 corresponding to 45 °. So on YBC 7289, we not only find the approximation for root 2 according to 1.24.51.10, but this is at the same time the function value of the secant $\alpha$ which by multiplying by the side length 30 and so in knowledge, but without using the Pythagorean theorem, into the real length 42.25.35 is converted \cite{fowler1998}.
\\
 
Based on our cuneiform tablet Plimpton 322, the gradient ratio 1 : 1 \ respective 45.00 ° would have been looked like this:

\begin{table}[h!]
\centering
\begin{tabular}{l l l l}
	\toprule
	\textbf{Column I} & \textbf{Column II} &  \textbf{Column III} & the unmentioned 3rd Side\\
	\midrule
	1.59.59.59.38.1.40 $\sim$ 2.0 & width = 6 & diagonal = 8.29.7 & (to the base = 6) \\
	
	59.59.59.38.1.40 $\sim$ 1.0 &  &  &\\
	\bottomrule
\end{tabular}
\end{table}

A shortening would also have been possible here, but the numerical values would not gave us a shorter written result. In addition to this mentioned last sentences, it becomes clear that Plimpton 322 only depicts a very limited and specially selected section of all valid gradient triangles. The reason for that will be possible explained at Point \{\ref{section:for?}\}.

\, \
\\

Another question that inevitably arises, especially in view of the attempts and looking for the solutions that have been going on for over 75 years is the following:

\subsubsection{Was the calculation of Pythagorean triples, or explicitly even their primitive variants the actual intention of the writers of Plimpton 322?}

No, because this is demonstrably not only a matter of the generally used rules, according to the Pythagorean theorem, but a much more complex variant of the equation derived from it. Likewise, it is not about the intention to generate Pythagorean triples in the table, the verifiable calculation step of expanding excludes this. Other Scientists in the field has mentioned this also, in her works about this cuneiform tablet\cite{Hajossy2016},\cite{abdulaziz2010},\cite{mansfield2017}.

It is to remind, that in principle every side of an exact right-angled triangle that can be given in sexagesimal numbers, by notation in decimal numbers is a triangle of Integer Sides, just depending on the starting Power to the initial sexagesimal place (1; 60; 3600; 216000). 
For that reason, even the Square root Values of Column I (as we know from chapter \{\ref{section:definitions}\}, these are function Values or ratios of $\dfrac{a}{b}$ otherwise $\dfrac{c}{b}$), are already and direct Pythagorean triangle sides in decimal notation. 

Conversely, this also means that, all sides in the relationship of a right-angled triangle to each other, in sexagesimal Numbers are also Pythagorean triangles. An additional extension by the reciprocal factor was not necessary for this.

For example at Line 1: The function Value $tangent$ or ratio $\dfrac{a}{b}$ is $59.30$ in decimal it is $3570$, and the Function Value $secant$ or ratio $\dfrac{c}{b}$ is $1.24.30$ in decimal $5070$. 
This is not a primitive however, a real Pythagorean Triangle with (Opposite, Hypotenuse, Adjacent) of  $(3570;5070;3600)$. The Babylonians has used instead the much more handy variant of $(119;169;120)$ by its smart extension. 
This chosen "smart" extensions made, by the stretching Factor $X$ the real side length $a$ for Column II, or $c$ for Column III from their ratios.
(see \{\ref{section:stepbystep}\})

Likewise, a purely mathematical / theoretical justification for the creation of the cuneiform tablet must be clearly denied. Although line 15 was unsuccessful due to the incorrect mathematical shortening, it is clear from both, line 11 and that line 15 that the author of the cuneiform tablet was not interested in the maximum possible mathematical shortening of his values. If possible, he also avoided a shortening of the base side of the right triangle below 1 = 60. This rather points to a geometrical, documentation and intended representative aspect, of the recorded triples. The absolute size of the resulting triangles was not decisive here, as long as the relationships were always maintained.

\, \

\subsubsection{Why do the gradient triangles contained on Plimpton 322 sometimes have larger gradations than the average? Was there perhaps no suitable triple for them?} 

This would speak against a systematic application of function values!

No, there are verifiable valid triples not only above and below the value range of line 1 or 15, but in addition in the larger gaps between (sometimes even several) also. Just for example, between lines 9 and 10, 11 and 12 or 14 and 15 \cite{kleb2020}. Whether, these were not selected from a point of view because of the given length in this case, or simply had no objective relevance for the content of this table, goes beyond the scope of this elaboration. In analyzing the column proportions of the tablet, however, one can assume that if it had been necessary, the Babylonian scribe could have noted also Numbers with a maximum of 10-11 places, instead of the noted 9 sexagesimal places in line 10, or the noted 8 sexagesimal places within line 2 (\{\ref{img001.jpg}\}).

\, \

\subsubsection{Was Plimpton 322 a school assignment, or an aid for a teacher to teach such incline relationships, as they are now demonstrably even contained in different forms or values on the cuneiform table?}

No, why do the angles only start at \textit{44.76 Degrees} and not with much more memorable gradient ratios and their current designation in degrees, for example also below \textit{53.13°}, which is a triple of $1.20; 1.40; (1)$ and would have been calculated according to the same rule\cite{Kleb2021a},\cite{Kleb2021b}. Especially for this mentioned incline relationship, there is evidence that many buildings were built not only in the Mesopotamian region after this\cite{gillings1982},\cite{bartclair2012}. Gradients higher than 1:1 there would only have been an automatic change of the short and longer triangle legs, but these could have been calculated just as naturally and in the same way as the existing value pairs and would have promoted understanding of the underlying Geometric.
 
Why should arithmetic problems for students only refer to closely staggered, steadily decreasing values and angles?(\cite{humml2019};\cite{businessinsider2018})

For learning about common relationships or their application in reality, or with common proportions and implicit angles of exactly \textit{30,45,60} or \textit{90} degrees or common inclines of 1:2 corresponding to \textit{26.56°} etc., no complex triples to be calculated were necessary. Or in some cases, small integer-numbered side lengths are not possible either. The easiest way to prove this, for example, is to use the vertical level with a plumb line and an isosceles triangle, which is used everywhere (builders level).

\, \

\subsubsection{Is the number of mathematical errors contained on the cuneiform tablet an indication for a work of a student?}

No, because as the analysis of the sources of error also presented in this paper clearly show \{\ref{section:errors}\} , a large part of the errors only arose during an actually unnecessary arithmetic operation. These were also calculated with a high degree of probability as mental arithmetic. Verifiable and in addition more difficult per line with a time offset between the calculation of the columns II and III. Perhaps the Babylonian scribe was under certain time pressure and maybe the lighting conditions were poor or unfavorable when writing. This is nothing unusual. It is concise, however, and so of special importance that they have a complex knowledge of the underlying geometry and mathematics in order to carry out the necessary operations. It corresponds to the Babylonian way of thinking, that controls were built in, for every value and the complete procedure\cite{Lemmermeyer2015a}\cite{lehmann1992}.

The author of the table knew this and regardless of column 0, which no longer exists today, the correct value of either other column could be reconstructed at any time and by everyone, with just a single value of the other columns. Therefore, he had to weigh up between calculation and writing speed with the possibility of later control and opposite the complete freedom from errors.

In this way, errors could have been corrected at the latest when applying or realigning the relevant gradient triangle, so that and this is a very important fact,"the value/gradient or Proportion definitely was not losing".

In this regard, I would like to point out that there is already clear evidence and a highly plausible theory, for the reasons of creating Plimpton 322 \{\ref{section:for?}\}, but in completeness this will be part of further research too.

\newpage
\section{Discussion of the Robson criteria}
\label{section:criteria}

In 2001, Eleanor Robson, as a known cuneiform-mathematics specialist set up six conditions or criteria\cite{robson2001}, according to which the respective statements or  new offered calculation procedures on the content of the cuneiform tablet should be checked; if these are not fulfilled, then the theory would exclude itself in particular for the considerations on trigonometric issues. So let's look at these conditions and look at their degree of fulfillment:

\subsubsection{1. The explanations must fit into the historical context of the ancient Babylonian period}

That the "Seked" described within the Papyrus Rhind (of the 16th century BC), as Cotangent, also a function value, like the discussed here, is an already accepted fact. Even that this Papyrus is also a copy of an older version from the 19th century BC \cite{mullerromer2016}\cite{bartclair2012}.

The usage of the reciprocal function value of $tangent$, is presented in detail on the basis of the works cited. 
A separation between the levels of knowledge of Egypt and Babylon, however, is not plausible.

The term "diagonal", which is often used in Old Babylonian, is a comparable term to describe a gradient ratio. In contrast to the defined height = $1$ at the Egyptian Seked in the gradient triangle\cite{gillings1982}\cite{furlong2019}. Against that, in old Babylonian the base was defined as $1$, which correlates as ratio with the tangent function value. The gradient ratio was, by that way, given as the mathematically unambiguous diagonal of this Gradient-Triangle.
From today's point of view, this apparently more abstract ratio of the $diagonals = secant$, compared to the $tangent$ as a direct gradient ratio and the $Seked = cotangent$ as its reciprocal, could represent the logical further development of the same functional value concept for extended applications.

To return to Babylon from Egypt, let us consider some recent scientific papers on astronomy and the use of mathematics and geometry in the late Babylonian period.
 
There are many mathematical cuneiform tablets which deal with the problems of inclines, proportions, cones, etc. Diagonals and a small Part outside the Triangle, but connected per definition, the Arrow, were explicitly used for this \cite{scriba2015,lehmann1992,robson2014}(for instance BM 085194,Db2-146). We find this especially together and connected with a geometry within a circle. 
On some of them, the ratios or proportions of the sides to one another, were multiplied by the real length to complete the task. Others, show calculations of uniform sectors in a circle and the use of circular arc segments. 

Furthermore, we know of extensive observations of astronomical events during this age \cite{vanderwaerden1942,weir1972,weir1982}, even if the way on collecting these observations has not yet been revealed, it is clear that there must have been very precise methods to locate special Locations at the Sky again and to compare with each other, even if there were years between. 

So it did not need our current angle naming practice to describe an incline or a proportionality to one another side length, for practical use.

The calculation sequence shown here is in that way completely in the historical context.

\subsubsection{2. The methodology presented must fit into the cultural consistency, here the general level of knowledge of the time and region is meant.}

In this case, the answer to the question overlaps with requirement or criteria 1. above. It can be stated, however, that the calculation sequence shown here is completely in the temporal and spatial context with various other time documents, but Plimpton 322 remains the most complex and therefore the unique of all these, worldwide.

\subsubsection{3. The presented arithmetic procedure, and other related analyzes, must fit into the arithmetic logic of the ancient Babylonian period and the region.}

Without repeating the calculations on Plimpton 322 again at this point, you can immediately see that some characteristics of the calculation can also be found on other cuneiform artifacts.(see also criteria 1. above and its references) Likewise, this procedure is also completely based on the specifications of the heading column I and is in its precision and abstraction ability, and only on the needs-based work, a very typical workflow for this time.

\subsubsection{4. The description according the presented procedure on Plimpton 322 must correspond to the physical conditions. This refers to the probably 5-6 cm larger horizontal dimension of the cuneiform tablet and that the method should therefore plausibly prove which data were on this lost part.}

As has been explained, the contents that were logically listed before today's column I fit in very well, both methodically and due to their brevity with a maximum of 4 sexagesimal places. In relation to the given proportions of today's remaining part and the possible original size, it can be said that this left lost part is slightly larger but at least as wide as today's columns 3 and 4. Based on the font size, it can even be anticipated that, there was another narrow Column, at the originally left edge of the tablet beneath column 0. This condition is therefore also to be regarded as fully met.

\subsubsection{5. Here the completeness and the direct reference between the headings and the numbers below, is postulated in a condition.}

The calculation procedure presented also completely fulfills this condition and shows us in a typical Babylonian way how intermediate calculations or further calculation steps should look like.

\subsubsection{6. The logical sequence of the table structure must correspond to the sequence of calculation steps and that there must not be any new content without reference to the other columns.}

This condition is also completely fulfilled based on the steps outlined.
\\

\textbf{Conclusion:} As the appreciated Eleanor Robson (\cite[199-202]{robson2001}) himself writes in her summary, should these conditions, be ``provocative'' also and additional aims into the direction of thinking of Buck, Friberg and Robson, this should exclude any legitimacy of trigonometric solutions. The level of fulfillment should be set so high that, these hurdles cannot be jumped over.

It was therefore at this time not actually intended that, the construction principle presented here fulfills all of these conditions, contrary to the assumptions made by the aforementioned 6 Criteria.
\\
\\

\section{What was the Plimpton 322 cuneiform tablet created for?}
\label{section:for?}

Since the system behind Plimpton 322 has now been extensively explained and it has also been shown that there was a continuous system of gradient ratios between 0 and 82 degrees, graduated in an average of 36 or less angular minutes, the primary question must be answered at this point. What was Plimpton 322 noted for and why does it contain these arbitrary angles? (also P. Moree 2018\cite{businessinsider2018}) 
On the other hand, it is clearly and not only immediately understandable for a surveyor that such nevertheless almost regular distances between the gradients / slopes of the documented triangles do not match and correlate to, terrestrial descriptions of a landscape, or a field. Even to buildings and city planning such regularities not really exist because of the always changing distances to the observed edges. 

Such terrestrial observed points on earth, could be only points that are far away on a circle and almost equally spaced. In addition, as an observer, you would have to have stand almost in the middle of this. However, the need and logic for such documentation seems to be completely incomprehensible. Therefore, we should look for other documented and plausible variants for such a procedure.
 
The only plausible variant, appears to be the documentation of a steady movement. Starting from an observation point, the position and thus also the height of an observed object changes in relation to a near fix time interval.  

The theory presented here is based on the knowledge of Babylonian mathematics and astronomy. We know that, at almost the same time, the Venus tablets of Ammi-Saduqa were created (Odenwald 2019\cite{odenwald2019} and others). The sidereal orbital period of the planet Venus is analyzed there by using of various observation cycles and this is set to 584 days. With today's measurement methods, the real value of such a cycle is 583.92 days \cite{vanderwaerden1942}. From this, we recognize the particularly high importance of astronomical observations in old Babylonian period. The close connection between astronomy and mathematics also becomes clear, but not only at that time millennia before us. In addition, international known Mathematicians of the 20th century have dealt intensively with those Venus tables. We can therefore be sure that even in ancient Babylonian times the best mathematicians, equipped with the best available possibilities and calculation methods, carried out the appropriate measurements. For such observations performed over many years with the highest precision, documentations such as the Plimpton 322 and high-precision time measurement was essential. Outside of daily use, water clocks were not useful for such complex observations.

\begin{figure}[h]
	\centering
	\includegraphics{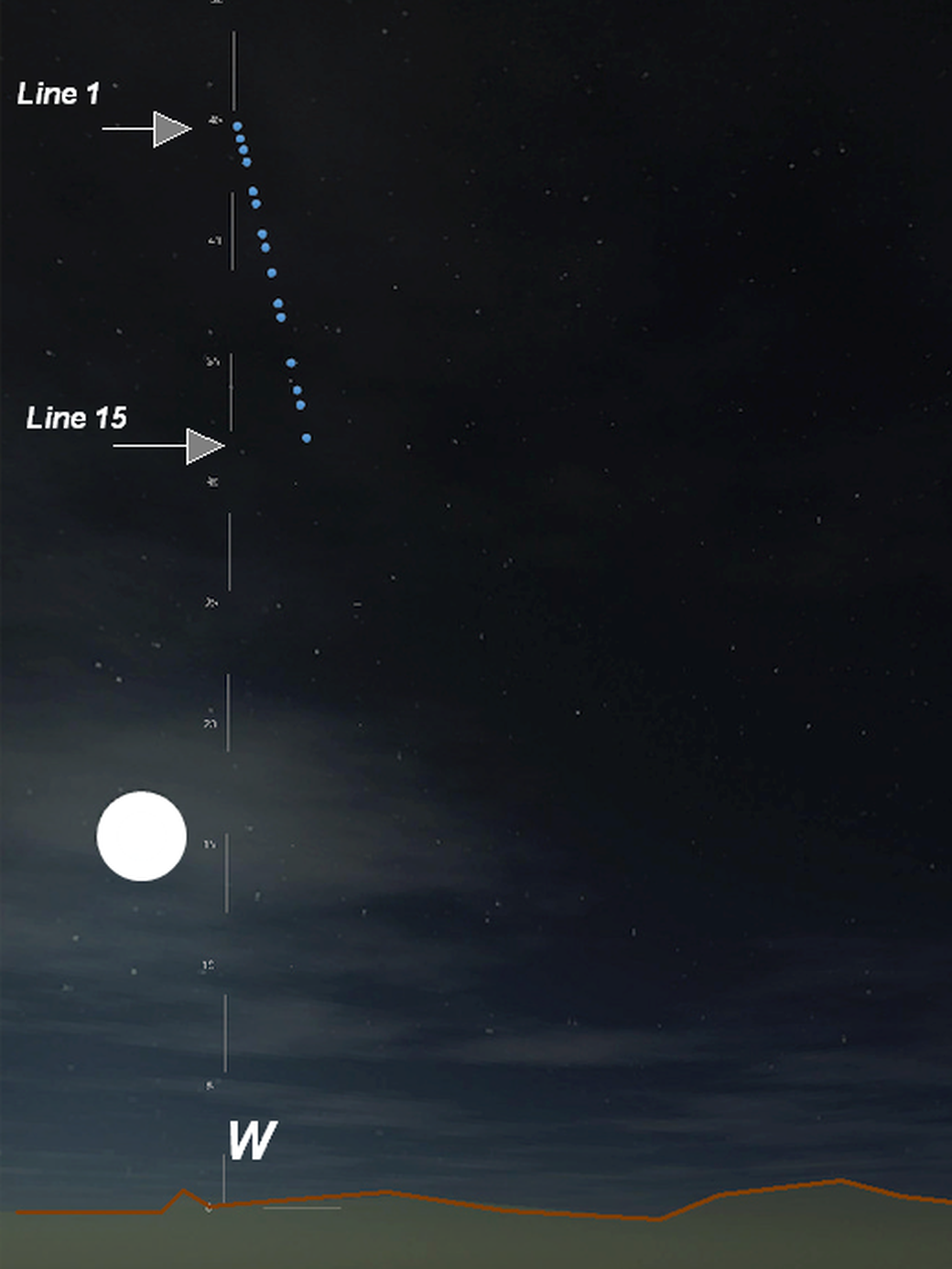}
	\caption{Given gradients of Line 1 to 15 embedded in a western evening Sky (unscaled drawing).}
	\label{fig:imagestars}
\end{figure}

Strongly connected to nature and deriving almost everything from it for practical everyday life, the Babylonians, to whom we owe the sexagesimal division of time and many more, demonstrably used astronomical occasions for the specification of time.

Grasshoff, Fatoohi and Stephenson, or Ossendrijver prove in their work that, at the latest from 700 BC, angle measurements with an accuracy of around 5-6 angular minutes were used for the documentation of astronomical phenomena \cite{fatoohi1998},\cite{grasshoff2017}. They pointed out that many of the cycles and observations described consisted of, or are derived from, much older data collections. M. Ossendrijver also describes exemplary but logically well-founded that there must have been appropriate measuring devices. That the astronomers, who were demonstrably servants or employees of the temple at this time, on the upper temple platform of the ziggurats, not only in the Neo-Assyrian and late Babylonian times, has carried out such observations appears just as plausible on closer inspection\cite{ossendrijver2013}\cite{ossendrijver2016}. In summary, it can be stated that it does not seem logical that the methodology of observation and measurement emerged spontaneously and without development. This is difficult to imagine and also not plausible in view of the parallel Egyptian knowledge.

HS0245 from the middle Babylonian period (1500-1000BC) presents us exact angular/visible distances between several stars given in values. These observations must have been carried out, as the stars and the moon, were high in the sky visible. This indicates even in the time before 700BC an observation equipment. (\cite{oracc2007},HS0245)

Since one can assume a steady refinement of the measuring methods and the used equipment over the centuries, it is therefore highly plausible to assume that, an older measuring method has already existed, before the known gradation of around 5-6 angular minutes (700 BC). 
With a much coarser gradation (by a factor of 6times coarser, corresponding to an average of 36 angular minutes), the calculation method presented here is a logical predecessor of the later ones, but it was already usable for the documentation of observations.

The Babylonian day began with the disappearance of the solar disk on the western horizon. The new month began, with the observation of exactly the same region of the sky, when the first delicate crescent of the moon short after new moon could be seen. Furthermore, as the sky disc of Nebra show us, at certain conjunctions between the Pleiades, the seven stars and the increasing moon, would again in the same celestial region and visible immediately after sunset observed that, an additional month have to be inserted, to mark the celestial passage of time and to guarantee its way.

We can thus state that the Babylonian time measurement was determined by astronomical events, as well as regularly adjusted by them. The same naturally applies to the astronomical measurements of that time already mentioned.

\begin{figure}[h]
	\centering
	\includegraphics[width=0.95\linewidth]{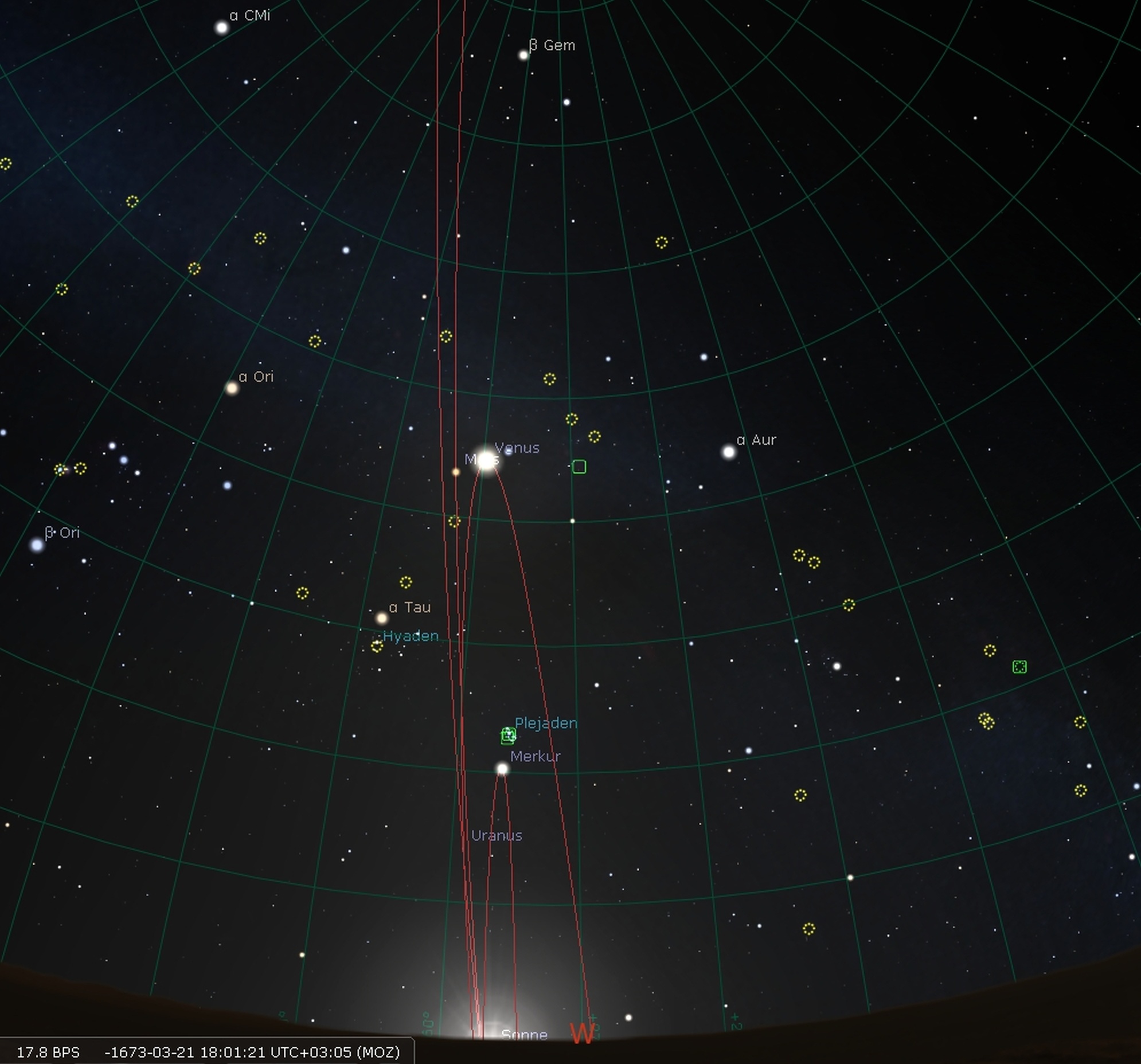}
	\caption{Calculated Venus-Position at Sunset, within the possible region (by Stellarium Software)}
	\label{fig:location1673bc}
\end{figure}

In fact, there is a regular astronomical event related to both, Venus and the precise line of the times, which was demonstrably of the greatest importance for the Babylonian astronomers. Here the planet Venus is at exactly the same position in the sky every 2919.6 days and once again this is the time of sunset, and at the same part in the western evening sky. More Precise, this event also took place in the vicinity of the Pleiades, which was already named with the Nebra Sky Disc.

If we do not attach any particular importance to the number 2919.6 days, those in the know and we immediately know after multiplying the correct days of our year by 8 that the two numbers are almost the same.

The 2919.6 Earth days correspond to 5 orbits of Venus, and 8 years of 365.26 days correspond to 2922.1 of our days. With a Precision after 8 years of 99.914\%, or a corresponding deviation of only 2.5 days (or 7.5 hours per year), this heavenly constellation regularly returned. With a deviation of around 3 seconds per hour from our normal time and guaranteed without consequential errors during the 8 years of operation, this event was certainly more precise than some new mechanical watches of our time.

This 8 year cycle, known to astronomers all over the world, always occurs at the time of the greatest eastern elongation of Venus and is one of the most striking astronomical events.

\begin{table}
	\centering
	\caption{Calculated positions of Sun, Venus and other Planets ( by Stellarium Software)}
	\label{venuspositions}

	\begin{tabular}{|l|l|l|l|l|l|l|l|}
	\hline  Name & Azimuth & Altitude & Mag. & Transit & Elev. at Transit & Elong. & Type \\ \hline
	\toprule
		\multicolumn{8}{l}{\cellcolor{midgrau}Positions at -1713-04-01 18:06:45 , moons age 29,1 days \,  \,  \,   \,   \   delta 1h 0m 11sec  }  \\ \hline
		Mercury & +268°37'49.8" & +21°01'42.7" & 0.93 & 13h22m & +323°26'10.53" & +21°27'14.3" & planet \\ \hline
		Sun & +268°24'59.8" & +0°04'49.5" & -21.90 & 12h09m & +45°35'19.71" & — & star \\ \hline
		Venus & +265°56'45.2" & +44°45'34.8" & -4.22 & 14h53m & +126°55'47.14" & +45°15'43.6" & planet \\ \hline
	\toprule	
		\multicolumn{8}{l}{\cellcolor{midgrau}Positions at -1705-03-30 18:05:49 , moons age 25,7 days  \,  \,  \,   \,   \  delta 1h 0m 12sec}  \\ \hline
		Mars & +262°47'01.4" & +37°43'01.9" & 1.93 & 14h31m & +180°02'56.18" & +38°25'28.0" & planet \\ \hline
		Sun & +267°32'33.5" & +0°04'18.8" & -21.88 & 12h09m & +2°23'39.59" & — & star \\ \hline
		Venus & +265°00'19.9" & +44°45'36.3" & -4.22 & 14h53m & +92°30'17.28" & +45°16'32.1" & planet \\ \hline
		\toprule
		\multicolumn{8}{l}{\cellcolor{midgrau}Positions at -1697-03-28 18:04:52 , moons age 22,2 days  \,  \,  \,   \,   \   delta 1h 0m 13sec}  \\ \hline
		Mars & +108°18'17.8" & +77°35'24.2" & 0.38 & 18h57m & +214°46'00.00" & +102°00'36.6" & planet \\ \hline
		Sun & +266°40'24.2" & +0°03'37.1" & -21.86 & 12h10m & +319°19'38.76" & — & star \\ \hline
		Venus & +264°02'41.1" & +44°45'31.4" & -4.23 & 14h54m & +57°27'32.50" & +45°17'31.9" & planet \\ \hline
	\toprule	
		\multicolumn{8}{l}{\cellcolor{midgrau}Positions at -1689-03-25 18:03:22 , moons age 17,7 days \,  \,  \,   \,   \ delta 1h 0m 19sec}  \\ \hline
		Sun & +265°21'17.3" & +0°03'19.5" & -21.85 & 12h11m & +254°27'43.58" & — & star \\ \hline
		Venus & +262°26'31.4" & +44°45'35.9" & -4.22 & 14h54m & +359°29'14.29" & +45°18'46.9" & planet \\ \hline
	\toprule	
		\multicolumn{8}{l}{\cellcolor{midgrau}Positions at -1681-03-23 18:02:23 , moons age 14,1 days  \,  \,  \,   \,   \  delta 1h 0m 21sec} \\ \hline
		Mercury & +263°55'28.7" & +6°53'56.0" & -0.65 & 12h37m & +41°19'35.72" & +7°18'09.1" & planet \\ \hline
		Moon & +81°15'36.5" & +8°06'46.2" & -11.69 & 0h04m & +80°37'10.17" & +171°50'09.4" & moon \\ \hline
		Sun & +264°29'23.2" & +0°02'51.7" & -21.84 & 12h12m & +211°45'45.16" & — & star \\ \hline
		Venus & +261°24'04.5" & +44°45'31.6" & -4.23 & 14h55m & +322°03'36.80" & +45°19'48.8" & planet \\ \hline
	\toprule	
		\multicolumn{8}{l}{\cellcolor{midgrau}\textbf{Positions at -1673-03-21 18:01:22} , moons age 10,7 days \,  \,  \,   \,   \   delta 1h 0m 27sec} \\ \hline
		Mars & +257°03'54.6" & +43°35'07.9" & 1.86 & 14h56m & +133°58'36.77" & +44°28'13.8" & planet \\ \hline
		Mercury & +263°45'07.6" & +20°20'24.3" & 0.73 & 13h24m & +77°57'22.01" & +20°48'34.9" & planet \\ \hline
		Moon & +99°04'33.9" & +48°22'14.2" & -11.43 & 21h01m & +42°42'05.08" & +130°19'27.1" & moon \\ \hline
		Sun & +263°37'54.7" & \textbf{+0°02'30.0"} & -21.82 & 12h13m & +169°29'02.58" & — & star \\ \hline
		Venus & +260°19'35.6" & \textbf{+44°45'35.3"} & -4.23 & 14h55m & +283°42'34.19" & +45°21'02.0" & planet \\ \hline
	\toprule	
		\multicolumn{8}{l}{\cellcolor{midgrau}Positions at -1665-03-18 17:59:41 , moons age 6,2 days \,  \,  \,   \,   \       delta 1h 0m 38sec} \\ \hline
		Mars & +92°51'39.9" & +65°16'29.5" & -0.02 & 19h51m & +231°58'06.67" & +114°47'11.2" & planet \\ \hline
		Moon & +225°53'50.9" & +71°35'53.7" & -9.72 & 17h03m & +61°41'37.01" & +75°46'24.2" & moon \\ \hline
		Sun & +262°19'13.5" & +0°02'47.9" & -21.84 & 12h14m & +105°23'37.17" & — & star \\ \hline
		Venus & +258°33'04.7" & +44°45'33.0" & -4.23 & 14h56m & +220°49'49.67" & +45°22'21.0" & planet \\ \hline
	\toprule	
		\multicolumn{8}{l}{\cellcolor{midgrau}Positions at -1657-03-16 17:58:37  ,  moons age 2,5 days  \,  \,  \,   \,   \   delta 1h 0m 50sec}  \\ \hline
		Moon & +262°08'03.9" & +30°03'02.5" & -6.63 & 13h49m & +293°24'41.40" & +30°31'55.1" & moon \\ \hline
		Sun & +261°28'32.5" & +0°02'44.9" & -21.84 & 12h15m & +63°59'30.17" & — & star \\ \hline
		Venus & +257°24'36.7" & +44°45'34.6" & -4.23 & 14h57m & +180°47'56.31" & +45°23'38.9" & planet \\ \hline
	\toprule	
		\multicolumn{8}{l}{\cellcolor{midgrau}Positions at -1649-03-14 17:57:27 , moons age 28,3 days \,  \,  \,   \,   \ delta 1h 0m 57sec} \\ \hline
		Sun & +260°37'30.5" & +0°02'43.7" & -21.84 & 12h15m & +22°20'43.90" & — & star \\ \hline
		Venus & +256°13'20.1" & +44°45'32.2" & -4.24 & 14h58m & +139°24'03.03" & +45°25'07.4" & planet \\ \hline
	\toprule
		\multicolumn{8}{p{16cm}}{\cellcolor{hellgelb} For Comparison some ``time to angular minutes'' Differences: 
			
			9sec. $\approx$ 2'  ;\,  \  1min. $\approx$ 13'  ;\,  \ 2m20sec. $\approx$ 30' = 0,5 Degrees
			
			delta = Time difference between the Gradients of Line 1 and Line 15 of Plimpton 322}
		\\\bottomrule
	\end{tabular}
\end{table}

If we now return to Plimpton 322 and the beginning of this chapter. The greatest eastern elongation of Venus, combined with the best and longest visibility, always occurs in March (beginning of April), i.e. near the vernal equinox. We remember the Babylonian beginning of the day, a moment in time as the last rays of the sun at the western evening sky were remains visible. In addition, this was a moment in time that, the Babylonians could estimate very precise. 

Corresponding to the former paragraphs we know that, any gradient triangle include an unique angle value. We also know that, there were valid gradient triangles for a wider range of angles, but at Plimpton 322 were used just a small range of them. The lines of decreasing gradient triangles at this cuneiform tablet starts with a triangle: $119; 169 ;(120)$ this is expressed in modern angles $\approx 44,7603 \text{° , which is } 44 \text{°}\, \ 45'\, \ 36"$. 

If we look at the day of the greatest eastern elongation of the Planet Venus exactly at sunset, so we will find the position exactly at an altitude of 44.76° above the horizon \{\ref{venuspositions}\}. This height above the horizon can be proved by software for historical astronomy, currently recommended by various observatories for the period of the 18th and 17th centuries BC. Due to the selection by the given calculation procedure, for the Babylonian gradient triangles and the so nearly fixed step size, we have to accept marginal inaccuracies of several time seconds or also angular minutes. 

In Addition, there are first answers about the initial intended observation interval. This ideal interval was 150 seconds, or 2 and half minutes. The real Observation for each Line, surely was within a little wider observation window, some seconds before and up to 10-15 seconds later, according to the used measuring procedure and equipment. 

The table \{\ref{venuspositions}\}, offers one but not the only one opportunity for the date of observation, it only should show the mentioned connection between sunset and the altitude of the Planet Venus at these events. Further comparisons with other cuneiform texts maybe could specify the data range much better.

Likewise, the classification of Plimpton 322 as documentation of a real series of observations leads to a further answer. The question was often asked why, because the angle differences were not regular and had \textit{gaps} in it. However, this is understandable for a \textit{real time} measurement which can not be repeated within a short period also if this is albeit an important one, because passing clouds or other events during the measurement could cause individual positions to \textit{fail}.

Due to the fixed gradation between the possible gradient ratios, caused to the selection of the values, a jump forward or backward would otherwise have had to take place in an intermediate measurement. But this would have led to a falsification of the measurement series in contrast to a simply larger interval.
 
Former NASA Scientist John D. Weir (\cite[30]{weir1982}) describes why the planetary movements in the western evening sky were used for such observations and that only the weather there could represent an obstacle to the detection of the position. Although the period for observation in the multi-year orbit cycle of the planet Venus, during the 2nd decade of March, is very special and unique, it can be stated the following. 
During this season of the year and especially in the region of the observation at that time, it is precisely provable that, this possible cloud cover during the observation can be plausibly explained. Based on modern weather observations in the region in question, we find 55-30\% cloud cover between February and April. (Source: weatherspark.com)

In this way, one of the last open questions about Plimpton 322 can be answered plausibly, the already pre-drawn lines, were on this way the preparation for further measurements up to the horizon. 

The fact that the series of values here, most probably due to the weather conditions has ended after the 15th observation (15 of initial intended 25, during the Observation-time between Line 1 and Line 15) and in addition at a very exact period of 1 hour (or 3600 seconds) \{\ref{venuspositions}\}, is certainly plausible for today's astronomers, or even understandable from their own practical experience. 

In summary, especially because that, this paper presents new answers, also further research will be necessary too. Maybe this Paper, will give by this way, some initial points for these.

\,
\,
\,

Germany, December 2021, Jens Kleb

\newpage
\nocite{*}
\begin{center}
	\textit{With many thanks, for the time-consuming support, to Patrick Leiverkus and my Wife Nicole.}
\end{center}
\,  \
\\
\,  \
\\
\,  \
\\
\,
\,
\,
\,
\begin{center}
Legend:	\{...\} refers to internal Chapters or tables, while [...] refers to the References below.
\end{center}
\printbibliography[heading=bibintoc]

\end{document}